\documentclass[11pt, a4paper]{article}
\usepackage{amsmath, amssymb, amsfonts, amsthm}
\usepackage[pdftex]{graphicx, color}

\newcommand{\EE}{{\bf{E}}}
\newcommand{\PP}{{\bf{P}}}
\newcommand{\Var}{{\bf{Var}}}

\newtheorem{tm}{Theorem}
\newtheorem{lemama}{Lemma}

\newtheorem{rem}{Remark}
\newtheorem{propopo}{Proposition}
\begin{document}

\parindent=0pt

\smallskip
\centerline{\LARGE \bfseries Local probabilities of randomly stopped sums}
\centerline{\LARGE \bfseries of power law lattice random variables}

\par\vskip 3.5em
\centerline{Mindaugas Bloznelis}
\vskip 0.5truecm
\centerline{mindaugas.bloznelis@mif.vu.lt}
\vskip 0.5truecm
\centerline{Institute of Computer Science,}
\centerline{Faculty of Mathematics and Informatics, Vilnius University}
\centerline{Naugarduko 24, Vilnius 03225, Lithuania}

\vskip 1truecm
\centerline{\it
Dedicated to Professor Vygantas Paulauskas on the occasion of his 75th birthday
}
\vskip 1truecm
\abstract{
 Let $X_1$ and $N\ge 0$ be integer valued 
power law random variables. 
For a randomly stopped sum $S_N=X_1+\cdots+X_N$ of 
independent and identically distributed copies of 
$X_1$ we establish a first order asymptotics of the local probabilities $\PP(S_N=t)$ as $t\to+\infty$.
Using this result we show $k^{-\delta}$, 
$0\le \delta\le 1$ scaling of the local clustering 
coefficient  
(of a randomly selected vertex of degree $k$)  
in a power law affiliation network.
}

\vglue2truecm
Keywords: Randomly stopped sum, local probabilities, power law, lattice random variables, clustering coefficient,
random intersection graph.

\vglue1truecm


\newpage

\section{Introduction}

\bigskip

Let $X_1,X_2,\dots$ be independent identically 
distributed random variables.
Let $N$ be a non-negative integer valued random variable 
independent 
of the sequence $\{X_i\}$. 
The randomly stopped sum $S_N=X_1+\dots+X_N$ 
is ubiquitous in many applications.  The tail probabilities 
$\PP(S_N>t)$ have attracted considerable attention in 
the literature
and their asymptotic behavior as $t\to+\infty$ 
 is quite well understood, 
see, e.g. \cite{AleskevicieneLeipusSiaulys},
\cite{DenisovFossKorshunov2010},
 \cite{EmbrechtsKM}, \cite{Foss} and references therein.
 Here we are interested in the asymptotic behavior of the {\it local 
 probabilities} $\PP(S_N=t)$.
 We assume that random variables $X_i$ are integer valued and $t$ is an integer.
Our study is motivated by
 several questions 
from the area of complex network modeling. An important class of complex networks  have (asymptotic) 
vertex
degree distributions
of the form $S_N$, where $X_i$ are integer valued and  $N$ and/or $X_i$ obey 
power laws. For this reason a rigorous 
analysis of 
network characteristics related to vertex degree
(clustering coefficients, degree-degree correlation)  
requires a good knowledge of the asymptotic behavior of the  
local 
probabilities  $\PP(S_N=t)$ as $t\to+\infty$,
\cite{BloznelisDegreeDegree},
\cite{BloznelisGJKR2015}, 
 \cite{BloznelisPetuchovas}. We will present  
applications in more detail after formulating our main results.
 
 In what follows we assume that for some $\alpha>1$
 we have
 \begin{equation}\label{V-04-11}
 \PP(X_1=t)=t^{-\alpha}L_1(t),
 \qquad
 t=1,2,\dots,
\end{equation}
where $L_1$ is  slowly varying at infinity. We will also be interested in the special case where
\begin{equation}
\label{XI-25-1}
L_1\bigl(tL_1^{1/(\alpha-1)}(t)\bigr)
\sim L_1(t).
\end{equation}
Here and below  
$f(t)\sim g(t)$ means  
$f(t)/g(t)\to 1$ as $t\to+\infty$. 
In the   particular case of (\ref{XI-25-1}),
where $L_1(t)\sim a$ as $t\to+\infty$ for some constant  $a>0$, we have 
\begin{equation}\label{V-04-11+a}
\PP(X_1=t)\sim a t^{-\alpha}.
\end{equation}

\noindent
For  $X_1$ having a finite first moment
we denote $\mu=\EE X_1$. 
In Theorem \ref{T1+} below 
we assume that
for some $\gamma>1$
\begin{equation}\label{VIII-30-1}
 \PP(N=t)= t^{-\gamma}L_2(t),
 \qquad
t=1,2,\dots,
\end{equation}
where $L_2$ is  slowly varying at infinity.

\begin{tm}\label{T1+}
Let $\alpha,\gamma>1$.
Assume that (\ref{V-04-11}) and  
(\ref{VIII-30-1}) hold and 
$\PP(X_1\ge 0)=1$.

\noindent 
(i) For $\gamma>\alpha$ and $\gamma>2$ we 
have 
\begin{equation}\label{V-04-12}
 \PP(S_N=t)\sim (\EE N) \PP(X_1=t).
\end{equation}

\noindent
(ii) For $\alpha>\gamma$ and $\alpha>2$ 
we have
\begin{equation}
\label{VIII-27-3}
\PP(S_N=t)\sim  \mu^{-1}
\PP\bigl(N=\lfloor t/\mu\rfloor\bigr).
\end{equation}
\noindent

\smallskip

\noindent
(iii) For $\alpha>2$, $\gamma\ge 2$ and $\EE N<\infty$ we have
\begin{equation}
\label{VIII-27-5}
\PP(S_N=t)
\sim
(\EE N)\PP(X_1=t)
+
\mu^{-1}\PP\bigl(N=\lfloor t/\mu\rfloor\bigr).
\end{equation}

\smallskip

\noindent
(iv) For $\alpha,\gamma<2$  conditions
(\ref{V-04-11+a}) and  (\ref{VIII-30-1}) imply
\begin{equation}\label{2017-XII-27-1}
\PP(S_N=t)
\sim 
t^{-1-(\alpha-1)(\gamma-1)} 
L_2(t^{\alpha-1})
a^{(\alpha-1)(\gamma-1)/\alpha}
(\alpha-1)
\EE Z_1^{(\alpha-1)(\gamma-1)}.
\end{equation}
Here  $Z_1$ is an $\alpha-1$ stable random variable with the characteristic function
\begin{displaymath}
\EE e^{i\lambda Z_1}
=
\exp
\left\{
|\lambda|^{\alpha-1}\Gamma(2-\alpha)
\left(\frac{\lambda}{|\lambda|}\sin\frac{(\alpha-1)\pi}{2}-\cos\frac{(\alpha-1)\pi}{2}\right)
\right\}.
\end{displaymath}
\end{tm}

In Theorems \ref{T2+} - \ref{RR4} below we 
drop the condition $\PP(X_1\ge 0)=1$. Instead we impose conditions on the left tail of $X_1$. 
In Theorem \ref{T2+} we focus on asymptotics (\ref{V-04-12}).
Assuming that (\ref{V-04-11}) holds we consider the following conditions on the left tail of $X_1$
\begin{eqnarray}
\label{2018-12-19-2}
&&
\PP(X_1<-t)=O\bigl(\PP(X_1>t)\bigr),
\\
\label{2019-01-06-1}
&&
\EE X_1^2{\mathbb I}_{\{X_1<-t\}} =o(\log^{-1}t)
\qquad
{\text{as}}
\qquad
t\to+\infty.
\end{eqnarray}
We also relax condition (\ref{VIII-30-1}) on the distribution of $N$ 
and consider the conditions 
\begin{eqnarray}\label{2017-06-04}
 &&
 \PP(N=t)=o\bigl(\PP(X_1=t)\bigr),
 \\
 \label{2019-01-07-1}
 &&
 \PP(N>t^2\ln^{-\tau}t)=o\bigl(t\PP(X_1=t)\bigr)
 \qquad
{\text{as}}
\qquad
t\to+\infty.
\end{eqnarray}

\begin{tm}\label{T2+}
Suppose that $\EE N<\infty$. 

\noindent 
(i) \hskip0.15truecm Let 
$1<\alpha<2$. Then  (\ref{V-04-11+a}) implies  (\ref{V-04-12}).

\noindent
(ii) \hskip0.05truecm
Let  $\alpha=2$. Assume that
$\EE (N\ln^{2+\tau}N)<\infty$, for some $\tau>0$. Then (\ref{V-04-11+a}) 
implies (\ref{V-04-12}).

\noindent
(iii) Let $2<\alpha<3$. For $\mu\le 0$ conditions (\ref{V-04-11+a}) 
and (\ref{2018-12-19-2}) imply
(\ref{V-04-12}). For $\mu >0$ conditions 

(\ref{V-04-11+a}), (\ref{2018-12-19-2})
and
(\ref{2017-06-04}) imply
(\ref{V-04-12}).

\noindent
(iv) Let  $\alpha=3$. Assume that (\ref{V-04-11+a}),
(\ref{2018-12-19-2}) hold. 
For $\mu>0$, respectively $\mu=0$, we assume  
in 

addition that 
$\PP(N=t)=o(t^{-3}(\ln\ln t)^{-1})$,
respectively  $\EE (N\ln^{1+\tau}N)<\infty$ 
for 
some 

$\tau>0$. Then (\ref{V-04-12}) holds.

\noindent
(v) Let  $\alpha>3$. Assume that (\ref{V-04-11}),
(\ref{2019-01-06-1}) hold. For $\mu=0$, respectively
$\mu>0$, we assume 
in

addition that
(\ref{2019-01-07-1}) holds
 for some 
$\tau>0$,
respectively (\ref{2017-06-04}) holds.
Then 
(\ref{V-04-12}) holds.

Moreover, for $\mu=0$ conditions
(\ref{V-04-11+a}),
(\ref{2019-01-06-1}) and (\ref{2019-01-07-1}) with $\tau=0$
imply 
(\ref{V-04-12}).
\end{tm}

\noindent
Theorem \ref{T2+} establishes (\ref{V-04-12}) under very mild conditions on $N$.
 For $\mu<0$  
 we only require the minimal condition
$\EE N<\infty$. 
Condition (\ref{2017-06-04}) and (\ref{2019-01-07-1}) with $\tau=0$ are minimal
ones 
as well.
The logarithmic factors 
in (\ref{2019-01-06-1}), (ii), (iv) are perhaps superfluous. 
They appear in the large deviation inequalities for
the tail probabilities $\PP(S_n>t)$ of sums $S_n=X_1+\dots+X_n$ that we apply in our proof.
Furthermore, condition (\ref{V-04-11+a}) on $X_1$ for 
$1<\alpha\le 3$ can be replaced by the weaker condition 
(\ref{V-04-11}), but then we need either an additional 
assumption 
on the slowly varying function $L_1$
or a bit stronger 
condition on $N$, see Theorems 
\ref{R2019-16-1}, \ref{R2019-16-2}  and
\ref{R1+NL1} below.
\begin{tm}\label{R2019-16-1} Let $1<\alpha\le 3$.
Assume that (\ref{V-04-11}) holds.
Assume that $\EE N^{1+\tau}<\infty$ for some $\tau>0$.
For $2<\alpha\le 3$ we also assume that 
(\ref{2018-12-19-2}) holds. For 
 $\mu>0$, $2<\alpha\le 3$ we assume in addition that
$\PP(N=t)=O(t^{-\alpha-\beta})$  for some $\beta>0$.
   Then  (\ref{V-04-12}) holds.
\end{tm}

\begin{tm}\label{R2019-16-2} 
Let $2< \alpha<3$. Suppose that $\mu>0$. 
Assume that  (\ref{V-04-11}), (\ref{XI-25-1}), 
(\ref{2018-12-19-2}), (\ref{2017-06-04})
hold. 
Then
 (\ref{V-04-12}) holds.
\end{tm}

In some applications a moment condition can be easier 
to verify than (\ref{2017-06-04}). In the following 
Theorem condition (\ref{2017-06-04}) is replaced  
by the moment condition $\EE N^{1+\alpha}<\infty$. Note 
that
(\ref{2017-06-04}) does not follow from  
$\EE N^{1+\alpha}<\infty$. 

\begin{tm}\label{R1+NL1}
Let $\alpha>1$. 
If $\EE N^{1+\alpha}<\infty$ then  (\ref{V-04-11+a}) implies (\ref{V-04-12}).
If $\EE N^{\beta}<\infty$ for some $\beta>1+\alpha$ then (\ref{V-04-11}) implies (\ref{V-04-12}).
\end{tm}

\noindent
It is interesting to compare the local probabilities of  $S_N$ 
with those of the maximal
summand $M_N=\max_{1\le i\le N}X_i$. Assuming that (\ref{V-04-11}) holds and $\EE N<\infty$
it is easy to show that
\begin{equation}\label{VIII-6-1}
 \PP(M_N=t)\sim (\EE N)\PP(X_1=t).
\end{equation}
Therefore, under conditions of Theorem \ref{T2+} the probabilities $\PP(S_N=t)$ and $\PP(M_N=t)$
are asymptotically equivalent.

Our next result Theorem \ref{T3+} establishes
(\ref{VIII-27-3}) assuming that $X_1$ is in the domain 
of attraction of a stable distribution with a finite 
first moment and 
\begin{equation}\label{VIII-27-2}
\PP(X_1=t)=o\bigl(\PP(N=t)\bigr)
\qquad
{\text{as}}
\qquad
 t\to+\infty.
\end{equation}
For $2<\alpha<3$ we will assume that 
\begin{equation}
\label{2018-12-19-1}
\exists 
\ 
\lim_{t\to+\infty}\frac{\PP(X_1<-t)}{\PP(X_1>t)}
<\infty.
\end{equation}
For $\alpha>3$ we assume (\ref{2019-01-06-1}). 
The case of $\alpha=3$ will be treated separately 
in Proposition \ref{X-12-1}.  
Assuming, in addition, that  $\EE N<\infty$ 
we show that as $t\to+\infty$
\begin{equation}\label{2019-02-11-1}
\PP(S_N=t)
=
(\EE N)\PP(X_1=t)(1+o(1))
+
\mu^{-1}
\PP\bigl(N=\lfloor t/\mu\rfloor\bigr)(1+o(1)).
\end{equation}
In our proof of  (\ref{VIII-27-3}), 
(\ref{2019-02-11-1}) 
we 
use 
regularity condition (\ref{VIII-30-1}).
In fact 
it can be slightly relaxed and replaced by the  following conditions
\begin{equation}\label{IX-11-1}
\exists 
\,
c_1>1,
\,
c_2,c_3>0:
\quad
 c_2\le \PP(N=t_2)/\PP(N=t_1)\le c_3
\quad
{\text{for any}}
\quad
1 \le t_2/t_1\le c_1,
\end{equation}
\begin{equation}\label{IX-11-2}
\exists
\,
\varkappa
=
\varkappa_{\alpha}
>
\max  \{(\alpha-1)^{-1}, 0.5\}:
\qquad
\lim_{t\to+\infty} \sup_{s:  |t-s|\le t^{\varkappa}}\PP(N=t)/\PP(N=s)=1.
\quad
\,
\
\
\end{equation}
In several cases we assume that for 
some $\beta>0$ 
\begin{equation}\label{2019-02-10-6}
\sum_{n\ge t}n^{-1-\beta}\PP(N=n)=o(\PP(N=t)).
\end{equation}
Clearly, (\ref{VIII-30-1}) implies 
(\ref{IX-11-1}),
(\ref{IX-11-2}), (\ref{2019-02-10-6})
but not vice versa.

\begin{tm}\label{T3+} 
Let $\alpha>2$. Suppose that $\mu>0$ and 
(\ref{IX-11-1}),
(\ref{IX-11-2}) hold.

\noindent
(i)
For $2<\alpha<3$ we assume that  (\ref{V-04-11+a}), 
(\ref{2018-12-19-1}) hold. Then  
$\EE N<\infty$ imply (\ref{2019-02-11-1}). 
Furthermore, if (\ref{V-04-11+a}), (\ref{VIII-27-2}),
(\ref{2018-12-19-1})  hold then 
either of the conditions $\EE N<\infty$  or 
(\ref{VIII-30-1}) with 
$\gamma>1$   imply
(\ref{VIII-27-3}).

\noindent
(ii)
For $\alpha>3$ we assume that  (\ref{V-04-11}),
(\ref{2019-01-06-1}) hold and  (\ref{2019-02-10-6})
is satisfied for 
$\beta=1/2$. 
Then  
$\EE N<\infty$ imply (\ref{2019-02-11-1}). 
Furthermore,
if in addition 
(\ref{VIII-27-2}) holds then 
either of the conditions $\EE N<\infty$  or 
(\ref{VIII-30-1}) with 
$\gamma>1$   imply
(\ref{VIII-27-3}).
\end{tm}

Theorem \ref{T3+} establishes (\ref{VIII-27-3}),
(\ref{2019-02-11-1}) under  
mild conditions on $N$.
  An inspection of the proof shows that
 (\ref{2019-02-10-6}) can be removed, but then we need  
 a much stronger condition (than (\ref{2019-01-06-1})) 
 on the rate of  decay of the left tail of $X_1$.
Condition 
(\ref{V-04-11+a}) 
of statement (i)
can also be slightly 
relaxed, see Theorem \ref{RR4} below.

\begin{tm}\label{RR4} Let $2<\alpha<3$. 
	 Suppose that $\mu>0$ and 
	(\ref{IX-11-1}),
	(\ref{IX-11-2}) hold.
Assume that (\ref{2019-02-10-6}) holds with some 
$\beta\in (0, \alpha-2)$. Then statement (i) of 
Theorem \ref{T3+}  remains valid if we replace 
condition (\ref{V-04-11+a}) by  
(\ref{V-04-11}), (\ref{XI-25-1}).
\end{tm}

\smallskip

\noindent
Now we consider the case where $\alpha=3$. We 
show in Proposition \ref{X-12-1} 
that
(\ref{V-04-12}), 
 (\ref{VIII-27-3})  and (\ref{VIII-27-5})  extend to 
 $\alpha=3$, 
 but we need a  
stronger condition on $X_1$ than (\ref{V-04-11+a}), (\ref{2018-12-19-1}). 
Namely, we assume that for some $a,b\ge 0$ such that $a+b>0$ and for some $\varepsilon>0$ we have for $t=1,2,\dots$
\begin{equation}\label{alpha=3}
\PP(X_1=t)
=
at^{-3}(1+r(t)),
\qquad
\PP(X_1=-t)=bt^{-3}(1+r(-t)),
\end{equation}
where
$r(s)
=
O\bigl( (\ln\ln |s|)^{-1-\varepsilon}\bigr)$ 
as $|s|\to+\infty$.
Note that the left relation of (\ref{alpha=3}) reduces to (\ref{V-04-11+a})  if we only require  $r(s)=o(1)$.

\begin{propopo}\label{X-12-1} 
Let $\alpha=3$. Assume that $\mu>0$ and (\ref{alpha=3}) holds.

\noindent
(i)
Suppose that 
(\ref{IX-11-1}),
(\ref{IX-11-2}) hold.
Then  
$\EE N<\infty$ imply (\ref{2019-02-11-1}). 
Furthermore, if (\ref{VIII-27-2}), (\ref{IX-11-1}),
(\ref{IX-11-2}) hold then 
either of the conditions $\EE N<\infty$  or 
(\ref{VIII-30-1}) with 
$\gamma>1$   imply
(\ref{VIII-27-3}).

\noindent
(ii) Suppose that $\EE N<\infty$. Then 
(\ref{2017-06-04}) implies (\ref{V-04-12}).
\end{propopo}

\smallskip

\noindent

Our proofs of Theorems \ref{T2+},\ref{R2019-16-1},\ref{R2019-16-2},\ref{T3+},\ref{RR4} and Proposition \ref{X-12-1}
combine the local limit theorem for the sums $S_n$
and  large
deviations inequalities for the tail probabilities 
$\PP(S_n>t)$. 
Note that   the values
$\alpha=2$ and $\alpha=3$ of the power law exponent (\ref{V-04-11}) are thresholds for the
centering and 
 scaling that make the sequence of distributions of 
$S_n$ 
tight.
Consequently,
large deviation inequalities at the threshold values
$\alpha=2$,
$\alpha=3$  have a 
slightly 
different form  comparing to the other (``ordinary'') values of 
$\alpha>1$. This is the (technical) 
reason  why 
for $\alpha=2$ and $\alpha=3$ we need more 
restrictive conditions on $X_1$ and/or 
$N$. In the proof of Theorem \ref{T1+} 
(where non-negative summands are considered)  
for $\alpha>2$  we use the discrete renewal theorem \cite{ErdosFellerPollard1949} and large deviation result \cite{Doney1989} for the local probabilities 
$\PP(S_n=t)$.

The results above are new. In the literature,
 see \cite{EmbrechtsKM},
\cite{Foss}, \cite{NgTang2004}, (\ref{V-04-12}) has been shown assuming a finite exponential moment $\EE e^{\delta N}<\infty$, for some 
$\delta>0$ (cf. Theorem  \ref{R1+NL1}, where the weaker moment condition  $\EE N^{1+\alpha}<\infty$   is assumed).
 (\ref{VIII-27-3}) has been established in \cite{Hilberdink1996} assuming that 
$\PP(X_1\ge 0)=1$ and
$\EE e^{\delta \sqrt{X_1}}<\infty$ for some $\delta>0$.
On the other hand
the asymptotics 
of the tail probabilities $\PP(S_N>t)$ 
corresponding to (\ref{V-04-12}), 
 (\ref{VIII-27-3})  and (\ref{VIII-27-5})
 in the presence of heavy tails are discussed in a number of papers, see \cite{AleskevicieneLeipusSiaulys},
\cite{DenisovFossKorshunov2010} and references therein.
We also mention a related work \cite{EmbrechtsMaejimaOmey1984} 
which establishes the limit
$\lim_{x\to+\infty}\PP(S_N\in (x,x+h])/\PP(N=\lfloor x\rfloor)$ in the case of non-negative 
  non-lattice summands  using
Blackwell's continuous renewal theorem. The  ratio
$\PP(S_N\in (x,x+h])/\PP(X_1\in (x,x+h])$  for $x\to+\infty$
is studied in  \cite{AsmussenFossKorshunov2003},
\cite{NgTang2004}, \cite{YuWangYang2010} 
in the case where $\EE e^{\delta N}<\infty$ for some 
$\delta>0$.


Before turning to applications we   mention two interesting  questions. The first question is about a higher order asymptotics, e.g., $k$ term   asymptotic expansion, to the probability $\PP(S_N=t)$ as 
$t\to+\infty$, for $k=2,3,\dots$. In the particular case of positive summands with $\EE e^{\delta \sqrt{X_1}}<\infty$ for some $\delta>0$ this problem is addressed in 
\cite{Hilberdink2009}, see also  \cite{BaltrunasSiaulys2007}.
 The second question is  what are the minimal conditions on the distributions of
  lattice random variables   $X_1$ and 
 $N$  satisfying (\ref{V-04-11}) and  (\ref{VIII-30-1}) that are sufficient 
 for either of 
 the relations (\ref{V-04-12}), 
 (\ref{VIII-27-3})  and (\ref{VIII-27-5}) to hold
 for various ranges of~$\alpha$.

\bigskip

{\bf Application to complex network modeling}. 
Mathematical modeling of complex networks aims at explaining and reproduction
of characteristic properties of large real world networks. We mention the power 
law degree distribution, short typical distances and clustering to name a few.
Here we focus on the clustering property meaning by this the tendency  of  nodes 
to cluster together by forming relatively small  groups with a high density of ties 
within a group. In particular, we are interested in the correlation between 
clustering and degree explained below.
Locally, in a vicinity of a vertex, clustering can be measured by the local clustering 
coefficient, the probability that two randomly selected neighbors of the vertex 
are adjacent. The average local clustering coefficient across vertices of degree 
$k$, denoted $C(k)$, for $k=2,3,\dots$, describes the correlation between 
clustering and degree.  Empirical studies of real social networks show that
the function $k\to C(k)$ is decreasing \cite{Foudalis2011}. Moreover, in the film 
actor network $C(k)$ obeys the scaling $k^{-1}$ \cite{RavaszB2003}. In the 
Internet graph it obeys the scaling $k^{-0.75}$ 
\cite{Vazquez2002}.  We are interested in modeling and explaining the scaling $k^{-\delta}$, for any given $\delta>0$ .
 
 \smallskip
 Clustering in a social network can be explained 
 by 
the auxiliary 
bipartite structure defining the adjacency relations between actors: every 
actor is  prescribed a collection of attributes and any two actors sharing an attribute have 
high chances of being adjacent, cf. \cite{NewmanStrogatzWatts2002}.  The 
respective random intersection graph  $G=G_{n,m}$ 
on the vertex set $V=\{v_1,\dots,v_n\}$ and 
with the auxiliary set of attributes $W=\{w_1,\dots, w_m\}$
 defines adjacency relations  between vertices with the help of a random 
bipartite graph $H$ 
linking actors to attributes. Actors/vertices are assigned iid non-negative weights 
$Y_1,\dots, Y_n$ modeling their activity and attributes are assigned iid 
non-negative weights $X_1,\dots, X_m$ modeling their attractiveness.
Given the weights, an attribute  $w_i$ is linked to actor $v_j$ in $H$ with 
probability $\min\{1, X_iY_j/\sqrt{mn}\}$ independently across the pairs $W\times V$. The 
pairs of vertices sharing a common neighbor in $H$ are declared adjacent in  
$G$.  The random intersection graph $G$ admits tunable power law degree 
distribution, non-vanishing global clustering coefficient, short typical distances, 
see  \cite{BloznelisGJKR2015}. 
Here we show that for large $m,n$  the random graph $G$ possesses yet another nice property, the  tunable scaling $k^{-\delta}$, $0\le \delta\le 1$, of  respective conditional probability 
$
C_G(k)= \PP
\bigl(
{\cal E}_{23}\bigr| {\cal E}_{12},{\cal E}_{13}, d_1=k
\bigr)$, a theoretical counterpart of $C(k)$.
By ${\cal E}_{ij}$ we denote the event that $v_i$ and $v_j$ are adjacent in $G$, $d_i$ stands for the degree 
of vertex $v_i$. In the following theorem, given two 
sequences of random weights 
$\{X_i, i\ge 1\}$ and $\{Y_j, j\ge 1\}$, we consider a family of random intersection graphs 
$\{G_{n,m}, n,m\ge 1\}$, where each $G_{n,m}$ is defined, by the weights $X_1,\dots, X_n$ and $Y_1,\dots, Y_m$ as above.

\begin{tm}\label{T4+}
Let 
$\alpha,\gamma>6$, $\beta>0$ and $a,\, b>0$. 
Suppose that $\{X_i, i\ge 1\}$ and $\{Y_j, j\ge 1\}$ are independent sequences of iid integer valued random variables  such that $\PP(X_i\ge 0)=1$, $\PP(Y_j\ge 0)=1$ and
$\PP(X_i=t)\sim a\,t^{-\alpha}$,  
$\PP(Y_j=t)\sim b\,t^{-\gamma}$ as $t\to+\infty$.
Let $m,n\to+\infty$. Assume that $m/n\to\beta$.
Then for every $k=2,3,\dots$ the probability
$C_G(k)
$ converges to a limit, denoted $C_*(k)$, and
\begin{equation}\label{X-24-1}
 C_*(k)\sim c k^{-\delta}
\qquad
{\text{as}}
\quad
k\to+\infty .
\end{equation}
Here
$\delta=\max\bigl\{0;\, \min\{\alpha-\gamma-1;1\}\bigr\}$,
$C_*(k)$ is given in  (\ref{X-28-2+}), and
$c>0$ is a constant depending on $\alpha,\gamma,\beta, a,b$ and the first three moments of $X_1$ and $Y_1$.
\end{tm}

A related result establishing $k^{-1}$ scaling in a random intersection graph with heavy tailed weights $Y_j$ and degenerate $X_i$ ($\PP(X_i=c)=1$ for some $c>0$) has been shown
in \cite{Bloznelis2013}. 
The {\it tunable} scaling $k^{-\delta}$, 
$\delta\in [0,1]$ in
(\ref{X-24-1}) is obtained due to the heavy tailed weights $X_i$. We suggest a simple 
explanation 
of how the weights of attributes  affect $C_*(k)$.  An attribute $w_i$ with weight $X_i$ generates with positive  probability
a clique in $G$ of size proportional to  $X_i$ (the clique formed by vertices linked to $w_i$).
For small $\alpha$ we will observe  quite a few large weights $X_i$. But the presence of many large cliques in $G$  may increase the value of $C_*(k)$ considerably. Hence, it seems plausible, that the scaling exponent $\delta$ correlated positively with $\alpha$.
For a  different approach to modeling of  
 $k^{-\delta}$ scaling, for $\delta=1$, we refer to  \cite{DorogovtsevGM2002}, \cite{RavaszB2003}.
 
 Another popular network characteristic that quantifies statistical dependence of neighboring adjacency relations is the correlation  coefficient (or rank correlation coefficient) between the degrees $d_1^*$ and $d_2^*$ of the endpoints of a randomly selected edge. More generally, one is interested in the distribution of the bivariate random vector $(d_1^*,d_2^*)$, called the "degree-degree" distribution. We briefly mention that using  the result of Theorem \ref{T1+} one  obtains from   Theorem 2 of \cite{BloznelisDegreeDegree}  that the  random intersection graph $G$ admits a tunable  power law degree-degree distribution. 
 
 The rest of the paper is organized as follows. Section 
 2 contains proofs of our main results. Large deviation 
 inequalities used in the proofs and an 
 upper bound for the convergence rate  in the local 
 limit theorem under condition (\ref{alpha=3}) are 
 given in Sections 3 and 4.

\section{Proofs}
Before the proofs we introduce some 
notation and present 
auxiliary lemmas. Then we prove our main results Theorems 
\ref{T1+}-\ref{RR4}
and Proposition
\ref{X-12-1}. 
Proofs of Theorem \ref{T4+} and relation (\ref{VIII-6-1}) are postponed to the end of this section.


{\bf Notation and auxiliary lemmas.} Given positive sequences $\{a_n\}$ and $\{b_n\}$ we denote $a_n\asymp b_n$ whenever $a_n=O(b_n)$ and $b_n=O(a_n)$ as $n\to+\infty$.
 We denote by 
 $c, c', c''$   positive constants, which may 
 depend on the distributions of $X_1$ 
 and $N$
 and 
may attain different values at 
different places.
 But they never depend on $t$.
For a non-random integer $n$ we 
 denote $S_n=X_1+\cdots+X_n$ and 
  $M_n=\max_{1\le i\le n}X_i$,
\begin{eqnarray}
\label{VII-11-5}
&&
\
S^{(1)}_n= \ X_1+\cdots+X_{\lfloor n/2\rfloor},
\qquad
S^{(2)}_n=X_{\lfloor n/2\rfloor +1}+\cdots+X_n,
\\
\nonumber
&&
M^{(1)}_n= \ \max_{1\le i\le \lfloor n/2\rfloor}X_i,
\qquad
\qquad
\
M^{(2)}_n= \ \max_{\lfloor n/2\rfloor<i\le n}X_i,
\\
\nonumber
&&
\,
Q^{(k)}_n= \ \sup_{i}\PP(S^{(k)}_n=i),
\qquad
\quad
L^{(k)}_n(t,\delta)=\PP\bigl(S^{(k)}_n\ge t/2,\, M^{(k)}_n< \delta t\bigr),
\
\
k=1,2.
\end{eqnarray}
Furthermore, in the case where $\EE X_1<\infty$, we 
denote  ${\tilde X}_i=X_i-\mu$ 
and ${\hat X}_i=\mu-X_i$.
We define ${\tilde S}_n$, 
${\tilde M}_n$, ${\tilde Q}^{(k)}_n$, ${\tilde L}^{(k)}_n$
in the same way as $S_n, M_n, Q^{(k)}_n, L^{(k)}_n$
above, 
but for the random variables ${\tilde X}_i$, $i\ge 1$. Similarly, 
we define ${\hat S}_n$,
${\hat S}^{(k)}_n$ in the same way as 
$S_n, S^{(k)}_n$ 
above,
but for the random variables ${\hat X}_i$, $i\ge 1$.
Given $t$ we 
denote 
$t_n=t-n\mu$ and ${\hat t}_n=n\mu-t$ so that
$\PP(S_n=t)=\PP({\tilde S}_n=t_n)=\PP({\hat S}_n={\hat t}_n)$.

In the proofs we bound the probability $\PP(S_n=t)$ by combining
two independent 
arguments: 
for large $n$ the probability is small by the local limit theorem and for large 
$t$ 
it is small because of the large deviations phenomenon. The argument is formalized in Lemma \ref{lema1+}.

\begin{lemama}\label{lema1+}  Let $0<\delta<1$. 
Let $n,t\ge 2$ be integers. 
We have
\begin{eqnarray}\label{lema1++}
 &&
\PP(S_n=t)
\le 
n\max_{i\ge \delta t}\PP(X_1=i)
+
\PP(S_n=t,\, M_n<\delta t),
\\
\label{lema1+++}
&&
\PP(S_n=t,\, M_n<\delta t)
\le
Q^{(1)}_nL^{(2)}_n(t,\delta)
+
Q^{(2)}_nL^{(1)}_n(t,\delta).
\end{eqnarray}
\end{lemama}
{\it Proof of Lemma \ref{lema1+}}. We have
\begin{displaymath}
 \PP(S_n=t)
=
\PP(S_n=t,\, M_n\ge \delta t)
+
\PP(S_n=t,\, M_n<\delta t).
\end{displaymath}
We evaluate the first probability on the right using the union bound, cf.  \cite{ZaigraevNJ},
\begin{eqnarray}\nonumber
\PP(S_n=t,\, M_n\ge \delta t)
&
\le
&
\sum_{1\le j\le n}\PP(S_n=t,\, X_j\ge \delta t)
=
n\PP(S_n=t,\, X_n\ge \delta t)
\\
\nonumber
&
=
&
n\sum_{i\ge \delta t}\PP(X_n=i)\PP(S_{n-1}=t-i)
\le
n\max_{i\ge \delta t}\PP(X_n=i).
\end{eqnarray}
It remains to  evaluate the second probability. We split
\begin{eqnarray}\nonumber
&&
 \PP(S_n=t,\, M_n<\delta t)
\\
\nonumber
&&
\qquad
\le
\PP(S_n=t, S^{(1)}_n\ge t/2, M_n<\delta t)
\ \,
+
\,
\PP(S_n=t, S^{(2)}_n\ge t/2, M_n<\delta t)
\\
\nonumber
&&
\qquad
\le
 \PP(S_n=t, S^{(1)}_n\ge t/2, M^{(1)}_n<\delta t)
+
\PP(S_n=t, S^{(2)}_n\ge t/2, M^{(2)}_n<\delta t)
\end{eqnarray}
and use the independence of 
$X_1,\dots, X_{\lfloor n/2\rfloor}$ and
$X_{\lfloor n/2\rfloor+1},\dots, X_n$. We have
\begin{eqnarray}\nonumber
 \PP(S_n=t, \, S^{(1)}_n\ge t/2, \, M^{(1)}_n<\delta t)
&
=
&
\sum_{i\ge t/2}\PP(S^{(2)}_n=t-i)\, \PP(S^{(1)}_n=i,\,  M^{(1)}_n<\delta t)
\qquad
\\
\nonumber
&
\le 
&
Q^{(2)}_n
\sum_{i\ge t/2}\PP(S^{(1)}_n=i, \, M^{(1)}_n<\delta t)
\\
\label{2019-01-14+1}
&
=
&
Q^{(2)}_nL^{(1)}_n(t,\delta).
\end{eqnarray}
We similarly show that
$
\PP(S_n=t, S^{(2)}_n\ge t/2, M^{(2)}_n<\delta t)
\le 
Q^{(1)}_nL^{(2)}_n(t,\delta)
$.$\qedsymbol$

 For integer valued iid random variables 
 $X_1,X_2,\dots$  satisfying
 (\ref{V-04-11}), (\ref{2018-12-19-1}) 
 the local limit theorem 
\cite{GnedenkoKolmogorov}, \cite{IbragimovLinnik}, \cite{Petrov} shows that
\begin{equation}\label{VIII-23-3}
\tau_n:=\sup_s\Bigl|b_n\PP(S_n=s)-g\bigl(b_n^{-1}(s-a_n)\bigr)\Bigr|\to 0
\qquad
{\text{as}}
\qquad
n\to+\infty.
\end{equation}
Here $b_n=n^{\beta}L_*(n)$ is a norming sequence,
$\beta=\max\{1/(\alpha-1); 0.5\}$
 and 
 $L_*(n)$ is a slowly varying function depending on 
 $\alpha$ 
and $L_1$, see (1.5.4), (1.5.5) in
 \cite{Borovkov2008}. 
For $\alpha<3$ we can choose
\begin{equation}\label{2018-I-2}
b_n=\max\bigl\{1, \inf\{x>0: \PP(|X_1|>x)<n^{-1}\}\bigr\}.
\end{equation}
For $\alpha>3$ we can choose $L_*(s)\equiv 1$. For 
$\alpha=3$ and $X_1$ satisfying (\ref{V-04-11+a}) we can choose 
$L_*(n)=\sqrt{\ln n}$.
Furthermore, $\{a_n\}$ is a centering sequence 
 ($a_n=0$ for $\alpha< 2$ and $a_n=n\mu$ for  $\EE |X_1|<\infty$), see, e.g., \cite{Borovkov2008}, and
$g(\cdot)$ is the probability density function of the stable limit 
distribution  of the sequence 
$\{(S_n-a_n)/b_n\}$.

\begin{lemama}\label{lema2++}  Let $\alpha>2$. 
Assume that (\ref{V-04-11}), (\ref{2018-12-19-1}) hold
and $\mu>0$.
 Then  as $t\to+\infty$
\begin{equation}\label{VIII-21-2}
\sum_{n:\, |n\mu-t|\le u_t}\PP(S_n=t)\to \mu^{-1}
\end{equation}
for any positive sequence $\{u_t\}$ satisfying 
\begin{equation}\label{VIII-23-4}
u_t/b_t\to +\infty,
\qquad
u_t^3b_t^{-2}t^{-1}\to 0,
\qquad
u_t\tau^*_{\lfloor t/(2\mu)\rfloor}/b_t\to 0. 
\end{equation}
Here $\tau_n^*:=\max\{\tau_k,\, k\ge n\}\to 0$ as $n\to+\infty$.
\end{lemama}

Notice that (\ref{VIII-23-4}) requires $u_t/b_t\to+\infty$ at a sufficiently slow rate.

{\it Proof of Lemma \ref{lema2++}}.
Denote for short ${\tilde t}=\lfloor t/\mu\rfloor$.
 Note that $a_n=n\mu$.

We establish (\ref{VIII-21-2}) in a few  steps
\begin{eqnarray}\label{VIII-23-5}
\sum_{n:\, |n\mu-t|\le u_t}\PP(S_n=t)
&
=
&
\sum_{n:\, |n\mu-t|\le u_t} b_n^{-1}g\bigl(b_n^{-1}(t-n\mu)\bigr)+o(1)
\\
\label{VIII-23-6}
&
=
&
\sum_{n:\, |n\mu-t|\le u_t} b_{n}^{-1}g\bigl(b_{\tilde t}^{-1}(t- n\mu)\bigr)+o(1)
\\
\label{VIII-23-8}
&
=
&
b_{\tilde t}^{-1} 
\sum_{n:\, |n\mu-t|\le u_t} g\bigl(b_{\tilde t}^{-1}(t- n\mu)\bigr)+o(1)
\\
\label{VIII-23-9}
&
=
&
\mu^{-1}+o(1).
\end{eqnarray}
Here (\ref{VIII-23-5}) follows from (\ref{VIII-23-3}) and the third relation of (\ref{VIII-23-4}). (\ref{VIII-23-6}) follows from the inequality shown below
\begin{equation}\label{VIII-23-10}
\Bigl|
\frac{1}{b_{\tilde t}}-\frac{1}{b_n}
\Bigr|
\le
c'\frac{u_t}{{\tilde t} b_{\tilde t}}
\end{equation}
combined with the mean value theorem (note that $g$ 
has a bounded derivative) and the second relation of (\ref{VIII-23-4}).
Furthermore, we obtain  (\ref{VIII-23-9}) by  approximating the sum by the integral of the unimodal 
density $g$ over the  unboundedly increasing 
domain $-u_tb_{\tilde t}^{-1}\le x\le u_tb_{\tilde t}^{-1}$
 \begin{equation}\label{VIII-23-11}
b_{\tilde t}^{-1}
\sum_{n:\, |n\mu-t|\le u_t}
 g\bigl(b_{\tilde t}^{-1}(t-n\mu)\bigr)
=
\mu^{-1}
\int_{|x|\le u_tb_{\tilde t}^{-1}}g(x)dx
+o(1)
\to
\mu^{-1}\int_{-\infty}^{+\infty}g(x)dx
=
\mu^{-1}.
\end{equation}
Finally,
(\ref{VIII-23-8}) follows from (\ref{VIII-23-10}) and (\ref{VIII-23-11}).

It remains to prove  (\ref{VIII-23-10}). We have
\begin{eqnarray}
\nonumber
&&
\frac{1}{b_{\tilde t}}-\frac{1}{b_n}
=
\Bigl(
\frac{1}{b_{\tilde t}}-\frac{1}{n^{\beta}L_*({\tilde t})}
\Bigr)
+ 
\Bigl(
\frac{1}{n^{\beta}L_*({\tilde t})}
-
\frac{1}{b_n}
\Bigr)
=:
I_1+I_2,
\\
\label{julius}
&&
|I_1|
=
\frac{|n^{\beta}-{\tilde t}^{\beta}|}
{n^{\beta}b_{\tilde t}}
\le
c'\frac{|n-{\tilde t}|\ {\tilde t}^{\beta-1}}
{n^{\beta}b_{\tilde t}}
\le
c'\frac{u_t}{{\tilde t} b_{\tilde t}},
\\
\label{VIII-23-2}
&&
|I_2|
=
\frac{1}{n^{\beta}}  \Bigl|\frac{1}{L_*({\tilde t})}-\frac{1}{L_*(n)}\Bigr|
\le
c'\frac{u_t}{b_{\tilde t}{\tilde t}}.
\end{eqnarray}
In (\ref{julius}) we applied the mean value theorem to 
$x\to x^{\beta}$.
In  (\ref{VIII-23-2}) 
we applied the inequality
\begin{equation}
|1-L_*(s+\delta_s)/L_*(s)|
\le
c'|\delta_s|s^{-1}
\end{equation}
to $s={\tilde t}$ and $s+\delta_s=n$.
To verify this inequality for large $s>0$ and $\delta_s=o(s)$ we use the 
representation $L_*(s)=c(s)e^{\int_{1}^s\varepsilon(y)y^{-1}dy}$, where
$\varepsilon(y)$ 
is a function satisfying $\varepsilon(y)\to 0$ as $y\to+\infty$,
 and where the  
$c(s)$
converges to a finite limit as $s\to+\infty$, see, e.g., \cite{Borovkov2008}.
Note that we can assume without loss of generality 
that $c(s)$ is a constant (as long as $L_*(n)$ defines a norming sequence). 
%
%
%
%
%
$\qedsymbol$

\begin{lemama}\label{clem} Let $\alpha>1$. 
Assume that  (\ref{V-04-11}) holds.
For any integers $n>0$ and $i$ we have  
$\PP(S_n=i)
\le c/b^+_n$. 
Here 
$\{b^+_n\}$  is the norming sequence of the sums 
$\{X_1^++\cdots+X_n^+, n\ge 1\}$ of iid random 
variables $X_1^+, X_2^+,\dots$ having the distribution
 $\PP(X_1^+=t)=\PP(X_1=t|X_1>0)$, $t=1,2\dots$, 
 cf. (\ref{VIII-23-3}), (\ref{2018-I-2}).
\end{lemama}
{\it Proof of Lemma \ref{clem}}.
Let $X_1^+$ and $X_1^{-}$ be  
random variables with the distributions 
\begin{displaymath}
\PP(X_1^+=t)=\PP(X_1=t|X_1>0), 
\qquad
\PP(X_1^-=-t)=\PP(X_1=-t|X_1\le 0),
\qquad
 t\ge 0.
\end{displaymath}
Let $\{X_i^+, i\ge 1\}$ and $\{X_i^{-}, i\ge 1\}$ be 
independent sequences of independent copies of $X_1^+$ 
and $X_1^{-}$.
Denote 
$n^{+}=\sum_{1\le i\le n}{\mathbb I}_{\{X_i>0\}}$.
For any integers $1\le r\le k\le n$ we have
\begin{equation}
\label{gruodzio18}
\PP(S_n=i|n^{+}=k)
=
\PP(X_1^{+}+\dots+X_k^{+}+X_{k+1}^{-}
+\cdots X_n^{-}=i)
\le 
Q_{r},
\end{equation}
where 
$Q_{r}
=
\sup_{j}
\PP(X_1^{+}+\dots+X_r^{+}=j)$. 
Furthermore, we have $\EE n^{+}=np$, where   
 $p=\PP(X_1>0)$. Put $r=\lfloor np/2\rfloor$. 
 Invoking Chernoff's bound 
$\PP(n^{+}\le r)\le e^{-r/4}$ and then
(\ref{gruodzio18}) we obtain 
\begin{equation}\label{gruodzio18+1}
\PP(S_n=i)
\le 
e^{-r/4}+\PP(S_n=i, n^{+}\ge r)
\le 
e^{-r/4}+Q_{r}.
\end{equation}
Note that for $\alpha\le 3$ ($\alpha>3$), $X_1^+$ is in 
the domain 
of attraction of $\alpha-1$ stable distribution (normal 
distribution). From  the local limit theorem we have 
$Q_r\le c/b_r^+$. Finally, the lemma follows from 
(\ref{gruodzio18+1}) combined with  relations  
$e^{-r}=o(1/b^{+}_r)$ and $b^{+}_r\asymp b^{+}_n$, for $r=\lfloor np/2\rfloor$.
$\qedsymbol$

\begin{lemama}\label{lema3}  Let $2\le \alpha<3$. Assume
 that $\EE|X_1|<\infty$ and $\mu>0$. Assume that (\ref{V-04-11}), (\ref{2018-12-19-1}), (\ref{XI-25-1}) hold and $\mu>0$. 
For $b_t$ defined by (\ref{2018-I-2}) and $A>1$ we have as $t\to+\infty$
\begin{eqnarray}\label{XI-25-3}
&&
\sum_{n:\, |n-t/\mu|\ge b_tA}|t_n|^{-\alpha}L_1(|t_n|)
\sim
c't^{-1}A^{1-\alpha},
\\
\label{XI-25-10}
&&
\sum_{n:\, |n-t/\mu|\ge b_tA}|t_n|^{-\alpha}\bigl(L_1(|t_n|)\bigr)^{\alpha/(\alpha-1)}
\sim
c't^{-1}A^{1-\alpha}L_*(t).
\end{eqnarray}
\end{lemama}
{\it Proof of Lemma \ref{lema3}}.
Recall that $b_t=t^{1/(\alpha-1)}L_*(t)$ and denote $t_{\star}=b_tA$.
 Note that (\ref{V-04-11}) implies
\begin{equation}\label{XII-05-1+}
\PP(X_1>t)
\sim 
t^{1-\alpha}
L_{1*}(t),
\qquad
{\text{where}}
\qquad
L_{1*}(t)
:=
(\alpha-1)^{-1}L_1(t).
\end{equation}
 Furthermore, (\ref{XI-25-1})
 implies
$L_{1*}\bigl(tL_{1*}^{1/(\alpha-1)}(t)\bigr)\sim L_{1*}(t)$,  because $L_1$ and $L_{1*}$ are slowly varying. Using Theorem 1.1.4 (v) of \cite {Borovkov2008}, we obtain from the latter relation that
\begin{equation}\label{L1*}
L_*(t)\sim L_{1*}^{1/(\alpha-1)}\bigl(t^{1/(\alpha-1)}\bigr).
\end{equation}
Recall that $t_n=t-\mu n$. Using  properties of slowly varying functions we evaluate the sums
\begin{eqnarray}\label{XI-25-4}
&&
\sum_{n:\, |n-t/\mu|\ge t_{\star}}|t_n|^{-\alpha}L_1(|t_n|)
\sim
c't_{\star}^{1-\alpha}
\
L_1(t_{\star}),
\\
\label{XI-25-6}
&&
\sum_{n:\, |n-t/\mu|\ge t_{\star}}|t_n|^{-\alpha}\bigl(L_1(|t_n|)\bigr)^{\alpha/(\alpha-1)}
\sim
c't_{\star}^{1-\alpha}
\
\bigl(L_1(t_{\star})\bigr)^{\alpha/(\alpha-1)}.
\end{eqnarray}
Furthermore, in view of (\ref{L1*}), we have
\begin{eqnarray}\nonumber
t_{\star}^{1-\alpha}
&
=
&
t^{-1}A^{1-\alpha}L_*^{1-\alpha}(t)
\sim 
c' t^{-1}A^{1-\alpha}L_1^{-1}\bigl(t^{1/(\alpha-1)}\bigr),
\\
\nonumber
L_1(t_{\star})
&
=
&
L_1\bigl(t^{1/(\alpha-1)}AL_*(t)\bigr) 
\sim
L_1\bigl(t^{1/(\alpha-1)}L_*(t)\bigr) 
\sim
L_1\Bigl(t^{1/(\alpha-1)}L_1^{1/(\alpha-1)}(t^{1/(\alpha-1)})\Bigr).
\end{eqnarray}
Combining  these relations with (\ref{XI-25-1}) we obtain
\begin{eqnarray}\label{XI-25-8}
&&
t_{\star}^{1-\alpha}
\ 
L_1(t_{\star})\sim  c't^{-1}A^{1-\alpha},
\\
\label{XI-25-9}
&&
L_1(t_{\star})
\sim
L_1\Bigl(t^{1/(\alpha-1)}L_1^{1/(\alpha-1)}(t^{1/(\alpha-1)})\Bigr)
\sim
L_1\bigl(t^{1/(\alpha-1)}\bigr)
\sim
c'L_*^{\alpha-1}(t).
\end{eqnarray}
In the very last step we applied (\ref{L1*}) once again.
Finally, invoking (\ref{XI-25-8}) in (\ref{XI-25-4}) and  (\ref{XI-25-8}), (\ref{XI-25-9}) in (\ref{XI-25-6}) we obtain (\ref{XI-25-3}), (\ref{XI-25-10}).
%
$\qedsymbol$

\noindent
{\bf Proofs of main results.} We first prove Theorems
\ref{T2+} - \ref{R1+NL1}
that establish (\ref{V-04-12}). 
The scheme of the proof is as follows. Given positive integer $m$  
 we split
\begin{eqnarray}
\label{b-17-1}
&&
\PP(S_N=t)
=
\EE \bigl(\PP(S_N=t|N)\bigr)
= 
I_m(t)+I'_m(t),
\\
\nonumber
&&
I_m(t)=\EE\bigl( \PP(S_N=t|N){\mathbb I}_{\{N\le m\}}\bigr),
\qquad
I'_m(t)=\EE\bigl( \PP(S_N=t|N){\mathbb I}_{\{N> m\}}\bigr)
\end{eqnarray} 
and show that  $I_m(t)=(1+o(1))(\EE N)\PP(X_1=t)$ and $I'_m(t)=o(\PP(X_1=t))$ for properly chosen  $m=m_t\to+\infty$ as $t\to+\infty$. We denote
\begin{displaymath}
  J_m=\EE\bigl(N{\mathbb I}_{\{N\le m\}}\bigr),
\qquad
J'_m=\EE\bigl(N{\mathbb I}_{\{N> m\}}\bigr).
\end{displaymath}
We begin with the proof of Theorem \ref{R1+NL1} since it is more transparent and simple.

{\it Proof of Theorem \ref{R1+NL1}}.
We recall the known fact that (\ref{V-04-11}) implies
for any $n=1,2,\dots$ 
\begin{equation}\label{V-04-10}
 \PP(S_n=t)\sim n\PP(X_1=t).
\end{equation}
Note that for any integer $m>0$, relation (\ref{V-04-10}) implies 
\begin{equation}\label{V-12-4}
I_m(t)\sim J_m \PP(X_1=t).
\end{equation}
It follows from   (\ref{b-17-1}) and (\ref{V-12-4}) that 
\begin{displaymath}
\liminf_{t\to+\infty}
\frac{\PP(S_N=t)}{\PP(X_1=t)}
\ge  J_m.
\end{displaymath}
Letting $m\to+\infty$ we obtain $J_m\to \EE N$ and 
\begin{equation}\label{V-12-1}
\liminf_{t\to+\infty}
\frac{\PP(S_N=t)}{\PP(X_1=t)}
\ge  \EE N.
\end{equation}
\noindent
To show the reverse inequality for $\limsup_t\bigl(\PP(S_N=t)/\PP(X_1=t)\bigr)$ 
we construct an upper bound for 
$I'_m(t)$. 

We first assume (\ref{V-04-11+a}) and $\EE N^{1+\alpha}<\infty$. 
By the union bound, 
\begin{eqnarray}\label{2019-01-07+2}
\PP(S_n=t)
&
\le
& 
n\PP(X_n\ge t/n, S_{n}=t)
=
n\sum_{i\ge t/n}\PP(X_n=i)\PP(S_{n-1}=t-i)
\\
\nonumber
&
\le
&
n\sup_{j\ge t/n}
\PP(X_n=j)\sum_{i\ge t/n}
\PP(S_{n-1}=t-i)
\le
n\sup_{j\ge t/n}\PP(X_n=j)
\
\
\
\\
\nonumber
&
\le
&
c'n(n/t)^{\alpha}.
\end{eqnarray}
Note that this inequality holds uniformly in $n$ and $t$.  Hence,
\begin{equation}\label{V-12-2}
I'_m(t)
\le 
c'
{\tilde J}_m'
t^{-\alpha},
\qquad
{\text{where}}
\qquad
{\tilde J}'_m=\EE\bigl(N^{1+\alpha}{\mathbb I}_{\{N> m\}}\bigr).
\end{equation}
It follows from (\ref{b-17-1}), (\ref{V-12-4}), (\ref{V-12-2}) that 
\begin{eqnarray}\label{VII-07+0}
\limsup_{t\to+\infty}\frac{\PP(S_N=t)}{\PP(X_1=t)}
&
\le
& 
\limsup_{t\to+\infty}\frac{I_m(t)}{\PP(X_1=t)}+\limsup_{t\to+\infty}\frac{I'_m(t)}{\PP(X_1=t)}
\\
\nonumber
&
\le
&
J_m
+
c'
{\tilde J}_m'.
\end{eqnarray}
Letting $m\to +\infty$ we obtain $J_m\to \EE N$ and ${\tilde J}_m'\to 0$.
Hence,
\begin{equation}\label{V-12-3}
\limsup_{t\to+\infty}\frac{\PP(S_N=t)}{\PP(X_1=t)}
\le 
\EE N.
\end{equation}
From (\ref{V-12-1}), (\ref{V-12-3}) 
we derive (\ref{V-04-12}).

Now we asume (\ref{V-04-11}) and $\EE N^{\beta}<\infty$, 
where $\beta>1+\alpha$. 
In view of (\ref{V-12-1}), (\ref{VII-07+0}) it suffices to show that
\begin{equation}\label{2019-10-11}
\lim_{m\to+\infty}
\limsup_{t\to+\infty}\frac{I'_m(t)}{\PP(X_1=t)}
=0.
\end{equation}
We write, for short,  $a_t^{1+\alpha}=t^{\alpha}L^{-1}_1(t)$ and   
split 
\begin{displaymath}
I_m'(t)
=
\EE
\left(
\PP(S_N=t|N)
\bigl(
{\mathbb I}_{\{m<N<a_t\}}
+
{\mathbb I}_{\{N\ge a_t\vee m\}}
\bigr)
\right)
=:J_{m.1}(t)+J_{m.2}(t).
\end{displaymath}
Next we prove  that 
$J_{m.1}(t)\le c\PP(X_1=t)\EE N^{\beta}{\mathbb I}_{\{m<N\}}$ 
and  
$J_{m.2}(t)\le \PP(X_1=t)\EE N^{1+\alpha}{\mathbb I}_{\{N\ge a_t\vee m\}}$. Note that these bounds imply (\ref{2019-10-11}).
We have, by Markov's inequality,
\begin{displaymath}
J_{m.2}(t)
\le 
\PP(N\ge a_t\vee m)
\le
a_t^{-1-\alpha}
\EE N^{1+\alpha}
{\mathbb I}_{\{N\ge a_t\vee m\}}
=
\PP(X_1=t)\EE N^{1+\alpha}{\mathbb I}_{\{N\ge a_t\vee m\}}
.
\end{displaymath}
To estimate $J_{m.1}(t)$ we use the property of a slowly varying function $L$ that for any $\delta>0$ there exists $y_0>0$ such that
$(x/y)^{\delta}\ge L(y)/L(x)\ge (y/x)^{\delta}$ for all $x\ge y\ge y_0$, see Thm. 1.1.2 and Thm. 1.1.4 (iii) 
in \cite{Borovkov2008}. For $\delta=\beta-1-\alpha$ and sufficiently large $t$ the left inequality implies 
$\max_{t/n\le j\le t}L_1(j)\le n^{\delta}L_1(t)$
for $n\le a_t=o(t)$. 
For $\delta=\alpha$ the right inequality  implies $\max_{j\ge t}j^{-\alpha}L_1(j)\le t^{-\alpha}L_1(t)$.
We have
\begin{equation}\nonumber
\max_{t/n\le  j\le t}\PP(X_1=j)
\le 
c(n/t)^{\alpha}
\max_{t/n\le j\le t}L_1(j)
\le 
cn^{\beta-1}t^{-\alpha}L_1(t)
\end{equation}
and
\begin{equation}\nonumber
\max_{ j\ge t/n}\PP(X_1=j)
\le 
\max_{t/n\le  j\le t}\PP(X_1=j)
+
\max_{ j> t}\PP(X_1=j)
\le
c(n^{\beta-1}+1)t^{-\alpha}L_1(t).
\end{equation}
Now (\ref{2019-01-07+2}) yields
$\PP(S_n=t)\le cn^{\beta}t^{-\alpha}L_1(t)$.
We obtain
$J_{m.1}(t)\le c\PP(X_1=t)\EE N^{\beta}{\mathbb I}_{\{m<N\}}$. 
$\qedsymbol$


{\it Proof of Theorem \ref{T2+}}.
We will construct a 
sequence $\psi_m=o(1)$ as $m\to+\infty$ and function $g(t)=o(\PP(X_1=t))$ such that
$I'_m(t)\le \psi_m\PP(X_1=t)+g(t)$ for $t>m$.
This inequality together with 
 (\ref{V-12-1}), (\ref{VII-07+0}) implies (\ref{V-04-12}).

We remark that (\ref{V-04-11}) implies 
$\PP(X_1=t)>0$ for sufficiently large 
$t$. For such $t$ we denote
$w_*(t)=\PP(N=t)/\PP(X_1=t)$ and 
$w(t)=\max\{w_*(s):\, s\ge t\}$.
Observe that (\ref{2017-06-04}) 
implies
$w(t)\downarrow 0$ as $t\to+\infty$. 
In the proof of 
(iii-v) below we assume that $t$ is sufficiently 
large so that $w(t)$  and
$w_*(t)$ 
are well defined. Denote $\delta=(\alpha-1)/(2\alpha)$.

{\it Proof of (i).} We estimate the probability $\PP(S_n=t)$ using
 Lemma \ref{lema1+}.
Invoking in (\ref{lema1++}) and (\ref{lema1+++})
the  inequalities  shown below
\begin{equation}\label{VII-1}
 \max_{j\ge \delta t}\PP(X_1=j)
\le
c't^{-\alpha},
\qquad
Q^{(k)}_n\le c'n^{-1/(\alpha-1)},
\qquad
L^{(k)}_n(t,\delta)
\le
c' n^{\alpha/(\alpha-1)}t^{-\alpha},
\end{equation}
we obtain $\PP(S_n=t)\le c'n t^{-\alpha}$. The latter inequality implies $I'_m(t)
\le  
c't^{-\alpha}
J'_m$. Clearly, $J'_m=o(1)$ as $m\to+\infty$.
It remains to prove (\ref{VII-1}). The first, second
 and third inequality of (\ref{VII-1})
 follows from (\ref{V-04-11}),
 Lemma \ref{clem} 
 and
(\ref{Borovkov1}) respectively.

{\it Proof of (ii).} Fix $\tau>0$. We show below that 
\begin{equation}\label{VII-3-2}
\PP(S_n=t)\le c't^{-2}n\ln^{2+\tau}n. 
\end{equation}
Note that (\ref{VII-3-2}) implies $I'_m(t)
\le  
c't^{-2}
\EE \bigl(N\ln^{2+\tau}N){\mathbb I}_{\{N\ge m\}}$
 and
 $\EE \bigl(N\ln^{2+\tau}N){\mathbb I}_{\{N\ge m\}}=o(1)$ 
as $m\to+\infty$.

\noindent
In the proof of (\ref{VII-3-2})
we 
distinguish two cases.
For $n\ln^{1+0.5\tau}n\ge t$ we have, by Lemma \ref{clem}, 
\begin{equation}\nonumber
 \PP(S_n=t)
\le
 c'n^{-1}
\le c'n^{-1} \frac{(n\ln^{1+0.5\tau}n)^2}{t^2}\le c't^{-2}n\ln^{2+\tau}n.
\end{equation}
For  $n\ln^{1+0.5\tau}n< t$ we show that 
$\PP(S_n=t)\le c'nt^{-2}$ using Lemma \ref{lema1+}
similarly as in the proof 
of statement (i) above. The only difference from (i) is that for $\alpha=2$ 
inequality (\ref{Borovkov1}), used in the proof of  
the third inequality of 
(\ref{VII-1}), only holds under additional condition (\ref{Borovkov1+}),
see Theorem \ref{Borovkovo} below. 
This condition (with $\beta<\tau/4$)
 is easily verified for  $x=t/2$, $y=t\delta$ 
and each $n$ satisfying
$n<t\ln^{-1-0.25\tau}t$. To derive the latter inequality from
$n\ln^{1+0.5\tau}n< t$ we argue by contradiction.
For $n_0\ge t\ln^{-1-0.25\tau}t$
we have (for sufficiently large $t$)
\begin{displaymath}
 n_0\ln^{1+0.5\tau}n_0
\ge
\frac{t}{\ln^{1+0.25\tau}t}\ln^{1+0.5\tau}\left(\frac{t}{\ln^{1+0.25\tau}t}\right)
=
(1+o(1)\bigr)t
\ln^{0.25\tau}t>t.
\end{displaymath}

\smallskip

{\it Proof of (iii)}. Assume that $\mu>0$.
We shall show that 
\begin{equation}\label{VII-11-1}
I'_m(t)
\le
c't^{-\alpha}
\Bigl(J'_m 
+ 
w^{1/2}\bigl(t/(2\mu)\bigr)\Bigr).
\end{equation}
We  
split
\begin{equation}\label{VI-29-1}
 I'_m(t)
=
\sum_{m< n<\infty}\PP( S_n=t)\PP(N=n)
=
I'_{m.0}+\dots+I'_{m.4},
\end{equation}
where $I'_{m.j}=I'_{m.j}(t)=\sum_{n\in{\cal N}_j}\PP(S_n=t)\PP(N=n)$ and where
\begin{eqnarray}
\label{VII-12-21}
&&
{\cal N}_0=\left( t_{-}; t_{+}\right),
\qquad
t_{\pm}=\frac{t}{\mu}
\pm t_{*},
\qquad 
t_*=
t^{1/(\alpha-1)}
w^{-1/2}\bigl(t/(2\mu)\bigr),
\\
\nonumber
&&
{\cal N}_1=\left(m; \frac{t}{2\mu}\right),
\quad
{\cal N}_2=\left[\frac{t}{2\mu}; t_{-}\right],
\quad
{\cal N}_3=\left[ t_{+}; \frac{2t}{\mu}\right],
\quad
{\cal N}_4=\left( \frac{2t}{\mu}; +\infty\right).
\end{eqnarray}
We obtain (\ref{VII-11-1}) from the bounds shown below
\begin{eqnarray}\nonumber
&&
I'_{m.0}
\le 
c'\PP(X_1=t)\, w^{1/2}\left(t/(2\mu)\right),
\\
\label{VII-11-3}
&&
I'_{m.j}
\le 
c't^{-\alpha}w^{(\alpha-1)/2}\left(t/(2\mu)\right),
\quad
j=2,3,
\qquad
I'_{m.j}\le c't^{-\alpha}J'_m,
\quad
j=1,4.
\qquad
\end{eqnarray}

\noindent
For $n\in{\cal N}_0$ we combine the bound 
$\PP({\tilde S}_n=t_n)\le c'n^{-1/(\alpha-1)}$ of Lemma \ref{clem} with (\ref{2017-06-04}) and obtain
\begin{eqnarray}\label{VII-30-3}
 I'_{m.0}
\le 
|{\cal N}_0|
\max_{n\in {\cal N}_0}\bigl\{\PP({\tilde S}_n=t_n)\PP(N=n)\bigr\}
\le 
c'\PP(X_1=t)\, w^{1/2}\left(\frac{t}{2\mu}\right).
\end{eqnarray}

Let us prove (\ref{VII-11-3}). 
We first show that 
\begin{eqnarray}\label{IX-08-1}
&&
\PP(S_n=t)
=
\PP({\tilde S}_n=t_n)\le c'nt_n^{-\alpha},
\qquad
{\text{for}}
\qquad
n\in {\cal N}_1\cup{\cal N}_2,
\\
\label{2019-01-10+1}
&&
\PP(S_n=t)
=
\PP({\hat S}_n={\hat t}_n)\le c'n{\hat t}_n^{-\alpha},
\qquad
{\text{for}}
\qquad
n\in {\cal N}_3\cup{\cal N}_4.
\end{eqnarray}
 We only prove (\ref{IX-08-1}). The proof of (\ref{2019-01-10+1}) is the same.
We apply Lemma \ref{lema1+} to the probability $\PP({\tilde S}_n=t_n)$. 
From   (\ref{lema1++}),
(\ref{lema1+++}) we obtain
\begin{eqnarray}\label{VII-11-2}
\PP({\tilde S}_n=t_n)
&
\le
&
n\max_{j\ge \delta t_n}\PP({\tilde X}_1=j)
+
{\tilde Q}^{(1)}_n{\tilde L}^{(2)}_n(t_n,\delta)
+
{\tilde Q}^{(2)}_n{\tilde L}^{(1)}_n(t_n,\delta)
\\
\nonumber
&
\le
&
 c'nt_n^{-\alpha}.
\end{eqnarray}
In the last step we used (\ref{VII-1}), which is shown 
using (\ref{Borovkov1}) similarly as in the proof of 
(i).

For $n\in {\cal N}_1$, 
respectively  $n\in {\cal N}_4$, 
the inequalities $t/2\le t_n<t$ and   (\ref{IX-08-1}),  
respectively
${\hat t}_n>t$ and (\ref{2019-01-10+1}), imply
\begin{equation}
\label{IX-17-1}
 \PP({\tilde S}_n=t_n)\le c'nt^{-\alpha},
\qquad
\PP({\hat S}_n={\hat t}_n)\le c'nt^{-\alpha}. 
 \end{equation}
  Hence the bounds 
$I'_{m.j}\le c't^{-\alpha}J'_m$, $j=1,4$.

For $n\in {\cal N}_j$, $j=2,3$ we use
$n\le c't$ and 
(\ref{IX-08-1}), (\ref{2019-01-10+1})
 to show that
\begin{equation}\label{IX-18-1++}
\sum_{n\in{\cal N}_j}\PP(S_n=t)
\le
c't\sum_{n\in{\cal N}_j}|t-n\mu|^{-\alpha}
\le
c't\,t_*^{1-\alpha}
\le c'w^{(\alpha-1)/2}\bigl(t/(2\mu)\bigr).
\end{equation}
Finally, the inequality
$\PP(N=n)\le c't^{-\alpha}$, 
$n\in{\cal N}_2\cup {\cal N}_3$, which follows from 
 (\ref{2017-06-04}), implies
\begin{equation}
I'_{m.j}
\le 
c't^{-\alpha}w^{(\alpha-1)/2}\bigl(t/(2\mu)\bigr),
\qquad
j=2,3.
\end{equation}
The proof for $\mu>0$ is complete.

Now assume that $\mu\le 0$.  
Inequalities $t_n\ge t$ and (\ref{VII-11-2}) imply 
$\PP({\tilde S}_n=t_n)
\le c'nt^{-\alpha}$, for $n\ge m$.
Hence $I'_m(t)\le c't^{-\alpha}J'_m$.
The proof for $\mu\le 0$ is complete.
\medskip

{\it Proof of (iv).} The proof is similar 
to that of (iii) above. One difference is that, 
for $\alpha=3$, 
large deviation inequality (\ref{Borovkov1}), 
used in (\ref{VII-11-2}), only holds under additional 
condition (\ref{2019-01-10+2}) that involves auxiliary 
functions $V(\cdot)$ and $W(\cdot)$, 
see Theorem \ref{Borovkovo} below.  
%
Let us now focus on (\ref{2019-01-10+2}). 
We note that (\ref{V-04-11+a}), (\ref{2018-12-19-2}) imply
\begin{equation}\label{2019-01-10+5}
\PP(X_1>u)\le c(1+u^2)^{-1},
\qquad
\PP(X_1<-u)\le c(1+u^2)^{-1},
\qquad u>0
\end{equation}
Using (\ref{2019-01-10+5}) we show that
$V(u)\le cu^{-2}\ln u$ and $W(u)\le cu^{-2}\ln u$
and reduce (\ref{2019-01-10+2})  to
\begin{equation}\label{2019-01-10+9}
 n
 \left(\frac{|\ln \Pi|}{y}\right)^2 
 \ln 
 \left(\frac{y}{|\ln\Pi|}\right)
\le \eta,
\qquad
{\text{where}}
\qquad
\Pi:=\Pi(x)=n(1+x^2)^{-1}.
\end{equation}

In the proof we also use the concentration bound
$\PP(S_n=j)\le c/\sqrt{n\ln n}$,
which follows from  (\ref{V-04-11+a}) by  
Lemma \ref{clem}. In particular, we have
${\tilde Q}^{(k)}_n$, 
${\hat Q}^{(k)}_{n}
 \le c/\sqrt{n\ln n}$, $k=1,2$.

\medskip

\noindent
Assume that $\mu<0$.
Using 
$t_n=t-n\mu\ge n|\mu|$ one easily shows that for some 
$\eta>0$ inequality
(\ref{2019-01-10+9}) holds with
$x=t_n/2$ and $y=t_n\delta$  uniformly in $n$. 
Now (\ref{Borovkov1}) implies 
\begin{equation}\label{2019-01-10+8}
{\tilde L}_n^{(k)}(t_n,\delta)
\le 
c'n^{\alpha/(\alpha-1)}t_n^{-\alpha},
\qquad
k=1,2.
\end{equation}
Invoking this bound in (\ref{VII-11-2}) and using
$t_n\ge t$ we obtain 
$\PP({\tilde S}_n=t_n)\le c'nt^{-\alpha}$.
Hence $I'_m(t)\le c't^{-\alpha}J'_m$. 
  The proof for $\mu<0$ is complete.

\medskip

\noindent
Assume that $\mu=0$. 
Denote $a_t=t^2\ln^{-1-0.5\tau}t$. We split
\begin{equation}\label{2019-09-16}
I'_m(t)
=
\Bigl(\sum_{m\le n< a_t}
+
\sum_{n\ge \max\{a_t,m\}}\Bigr)\PP(S_n=t)\PP(N=n)=:J_1+J_2
\end{equation}
and estimate
\begin{eqnarray}\label{2019-01-10+6}
&&
J_1
\le 
\sum_{m\le n< a_t}
c'nt^{-\alpha}\PP(N=n)
\le 
c'
J'_mt^{-\alpha},
\\
\label{2019-01-10+7}
&&
J_2
\le 
\frac{c}{\sqrt{a_t\ln a_t}}
\PP(N\ge a_t)
\le 
\frac{c}{a_t^{3/2}\ln ^{1.5+\tau}a_t}
\EE N(\ln N)^{1+\tau}{\mathbb I}_{\{N\ge a_t\}}
=
\frac{o(1)}{t^3\ln^{\tau/4}t}.
\quad
\end{eqnarray}
We obtain
 $I'_m(t)\le c'J'_mt^{-3}
 +
 o\bigl(t^{-3}\ln^{-\tau/4}t\bigr).
 $

 \noindent
 It remains to prove  (\ref{2019-01-10+6}), (\ref{2019-01-10+7}).
 In  (\ref{2019-01-10+7})
we first 
estimated  
\begin{displaymath}
\PP(S_n=t)\le c(n\ln n)^{-1/2}\le c
(a_t\ln a_t)^{-1/2},
\qquad
n\ge a_t,
\end{displaymath}
and then applied Markov's inequality. In (\ref{2019-01-10+6}) 
we invoked 
(\ref{VII-11-2}) with $t_n=t-n\mu=t$. Note that (\ref{VII-11-2}) uses the bound 
(\ref{2019-01-10+8}),
which follows from (\ref{Borovkov1}). But (\ref{Borovkov1}) only holds if condition (\ref{2019-01-10+9}) is satisfied.
Now we show that for some $\eta>0$ 
inequality 
(\ref{2019-01-10+9}) holds with $x=t/2$, $y=t\delta$ uniformly in
$n,t$ satisfying $n<a_t$. 
We consider the cases $n\le b_t:=t^2\ln^{-3}t$  and 
$b_t\le n\le a_t$ separately.
For $n\le b_t$ we have $c't^{-2}\le \Pi\le c''\ln ^{-3}t$ so that
the left side of (\ref{2019-01-10+9}) does not exceed
$c'n
 (t^{-1}\ln t)^2\ln t\le c''$.
For $b_t<n\le a_t$ we have  $c'\ln^{-3}t\le \Pi\le c''\ln^{-1-0.5\tau}t$ so that  the left side of (\ref{2019-01-10+9}) does not exceed
$c'n
 (t^{-1}\ln\ln t)^2\ln t\le c''$.
The proof for $\mu=0$ is complete.

\medskip

\noindent
Assume that $\mu>0$. We have that
$\psi_*(r):
=
\max_{t\ge r}\bigl(\PP(N=t)t^3\ln\ln t\bigr)\downarrow 0$
as $r\to+\infty$.
We shall show that 
\begin{equation}\label{2019-01-11+1}
I'_m(t)
\le 
c't^{-\alpha}
\bigl(J'_m+\psi_*(t/(2\mu))+w(t/(2\mu))\bigr).
\end{equation}
We put $t_*=\sqrt{t\ln t}\,(\ln\ln t)$ in (\ref{VII-12-21}), decompose 
$I'_m(t)$ using (\ref{VI-29-1}) and  obtain (\ref{2019-01-11+1})  from  the bounds shown below
\begin{eqnarray}\label{2019-01-11+3}
&&
I'_{m.0}
\le 
c't^{-\alpha}\psi_*(t_{-}),
\\
\label{2019-01-11+2}
&&
I'_{m.j}
\le 
c't^{-\alpha}w(t/(2\mu)),
\quad
j=2,3,
\qquad
I'_{m.j}\le c't^{-\alpha}J'_m,
\quad
j=1,4.
\qquad
\end{eqnarray}
We obtain (\ref{2019-01-11+3}) proceeding as in (\ref{VII-30-3}) and using 
$\PP(N=t)\le \psi_*(t_-)/(t^3\ln\ln t)$, 
$n\in{\cal N}_0$. To show (\ref{2019-01-11+2}) for $j=2,3$ we note that $\PP(N=t)=o(\PP(X_1=t))$ and combine the inequalities
\begin{displaymath}
\max_{n\in{\cal N}_2\cup{\cal N}_3}\PP(N=n)
\le 
w(t/(2\mu))
\max_{n\in{\cal N}_2\cup{\cal N}_3}\PP(X_1=n)
\le c w(t/(2\mu)) t^{-3}
\end{displaymath}
with the bound
$\sum_{n\in{\cal N}_j}\PP(S_n=t)
\le 
c't t_*^{1-\alpha}$ that follows from  
(\ref{IX-08-1}), (\ref{2019-01-10+1}), 
see  (\ref{IX-18-1++}). 
For $j=1,4$  the bounds (\ref{2019-01-11+2}) are derived from (\ref{IX-08-1}), (\ref{2019-01-10+1}) in the same way as those of (\ref{VII-11-3}).

Inequalities (\ref{IX-08-1}), (\ref{2019-01-10+1})
used in the proof of (\ref{2019-01-11+2}) above are 
derived from (\ref{Borovkov1}). In the case of   
(\ref{IX-08-1}), respectively (\ref{2019-01-10+1}),
we apply (\ref{Borovkov1}) to 
${\tilde S}_n$, $x=t_n/2$, $y=t_n\delta$, 
$n\in {\cal N}_1\cup{\cal N}_2$, respectively 
${\hat S}_n$, $x={\hat t}_n/2$, 
$y={\hat t}_n\delta$, 
$n\in {\cal N}_3\cup{\cal N}_4$.
But for $\alpha=3$ (\ref{Borovkov1}) only holds when 
condition (\ref{2019-01-10+2}) is satisfied.

In order to verify (\ref{2019-01-10+2}) in the case of 
(\ref{IX-08-1})
we show that for some $\eta>0$ inequality
(\ref{2019-01-10+9}) is satified by  $x=t_n/2$, $y=t_n\delta$, $n\in {\cal N}_1\cup{\cal N}_2$. Denote
$a_t=\lfloor t\mu^{-1}-t_*\ln t\rfloor$.
We consider the cases 
$n\in (m;a_t]$ 
and 
$n\in (a_t; t_{-}]$ 
separately. For
$n\in (a_t; t_{-}]$ 
we have
$c'\ln\ln t\le |\ln \Pi| \le c''\ln\ln t$.
Hence the left side of (\ref{2019-01-10+9}) 
is at most $n(t_n^{-1}\ln\ln t)^2\ln t
\le 
t_{-}((\mu t_*)^{-1}\ln\ln t)^2 \ln t\le c$.
For $n\in (m;a_t]$
 we have 
$c'\ln\ln t\le |\Pi|\le c''\ln t$.
Hence the left side of (\ref{2019-01-10+9}) 
is at most 
$n(t_n^{-1}\ln t)^2\ln t
\le 
a_t(t_{a_t}^{-1}\ln t)^2 \ln t\le c$.
Here we used the fact that $n\to nt_n^{-2}$ is increasing.

In order to verify (\ref{2019-01-10+2}) in the case of 
(\ref{2019-01-10+1})
we show that for some $\eta>0$ inequality
(\ref{2019-01-10+9}) with $\Pi=n\PP({\hat X}_1>x)$
is satisfied by 
$x={\hat t}_n/2$, 
$y={\hat t}_n\delta$, 
$n\in {\cal N}_3\cup{\cal N}_4$.
 Denote
$b_t=\lfloor t\mu^{-1}+t_*\ln t\rfloor$, 
$e_t=\lfloor 2t/\mu\rfloor$.
We consider the cases 
$n\in (t_{+};b_t]$, 
$n\in (b_t; e_t)$ 
and 
$n\ge e_t$
separately. For
$n\in (t_{+}; b_t]$ 
we have
$c'\ln\ln t\le |\ln \Pi| \le c''\ln\ln t$.
Hence the left side of (\ref{2019-01-10+9}) 
is at most $n({\hat t}_n^{-1}\ln\ln t)^2\ln t
\le 
b_t((\mu t_*)^{-1}\ln\ln t)^2 \ln t\le c$.
For $n\in (b_t;e_t)$
 we have 
$c'\ln\ln t\le |\Pi|\le c''\ln t$.
Hence the left side of (\ref{2019-01-10+9}) 
is at most 
$n({\hat t}_n^{-1}\ln t)^2\ln t
\le 
e_t({\hat t}_{b_t}^{-1}\ln t)^2 \ln t\le c$.
Finally, for $n\ge e_t$ we have 
$n\mu/2\le {\hat t}_n\le n\mu$ and
$c'\le n\Pi \le c''$.
Hence the left side of (\ref{2019-01-10+9}) 
is $O(n^{-1}\ln^3n)<c$.
The proof for $\mu>0$ is complete.

\smallskip

{\it Proof of (v).} 
Denote $\sigma^2=\Var X_1$,
$S_n^*=S_n-\lfloor n\mu\rfloor$, 
$t_n^*=t-\lfloor n\mu\rfloor$ so that 
$\PP(S_n=t)=\PP(S_n^*=t_n^*)$.

Assume that $\mu<0$. 
For $n>\mu^{-2}$
we have $t_n^*=t+|\lfloor n\mu\rfloor|\ge \sqrt n$ and, 
by (\ref{Borovkov8}),
\begin{equation}\label{2019-01-12+1}
\PP(S_n=t)
=
\PP(S_n^*=t_n^*)
\le
cn^{-1/2}e^{-(t_n^*)^2/(2n\sigma^2)}
+
cn(t_n^*)^{-\alpha}L_1(t_n^*).
\end{equation}
Here we assume that $n$ is sufficiently large so that
$(t^*_n)^{-1}\PP(X_1-\mu>t^*_n)
\le 
c (t^*_n)^{-\alpha}L_1(t^*_n)$,
see Thm.1.1.4 (iv) of \cite{Borovkov2008}.
Using 
$(t_n^*)^2\ge 2t|\lfloor n\mu\rfloor|
\ge 2tn|\mu|$ we estimate the first term~on~the~right
\begin{displaymath}
e^{-(t_n^*)^2/(2n\sigma^2)}
\le 
e^{-t|\mu|/\sigma^2}
\le 
c
t^{-\alpha}L_1(t).
\end{displaymath}
In the last step we assumed that $t$ is sufficiently large.
Furthermore, it follows from general properties of 
slowly varying functions that for some $s_0>0$ and 
$c>0$ (both depending on $L_1$) we have
$x^{-\alpha}L_1(x)\le c y^{-\alpha}L_1(y)$ for
 $x>y>s_0$. Hence, for sufficiently large $t$ we have
$(t_n^*)^{-\alpha}L_1(t_n^*)\le ct^{-\alpha}L_1(t)$ 
for each $n$.
Now (\ref{2019-01-12+1}) implies 
 $\PP(S_n=t)\le c'nt^{-\alpha}L_1(t)$. From the latter inequality  we obtain (for large $m$) 
$I'_m(t)
\le  
c't^{-\alpha}
J'_m$.
The proof for $\mu<0$ is complete.

 Assume that $\mu=0$ and (\ref{V-04-11}) holds. 
 Denote $b_t=t^2/(1+\sigma^2)$. For $n\in (m;b_t)$
 we estimate $\PP(S_n=t)$ using (\ref{Borovkov8}). For $n\ge b_t$ we apply
  the concentration bound $\PP(S_n=t)\le cn^{-1/2}$.
We obtain
\begin{eqnarray}\nonumber
&&
I'_m(t)
\le 
c(I_1+I_2+I_3),
\\
\nonumber
&&
I_1=\sum_{n\in (m;b_t)}n^{-1/2}e^{-t^2/(2n\sigma^2)}\PP(N=n),
\qquad
I_2=\sum_{n\in (m;b_t)}nt^{-1}\PP(X_1>t)\PP(N=n),
\\
\nonumber
&&
I_3=\sum_{n\ge b_t}n^{-1/2}\PP(N=n)
\le 
b_t^{-1/2}
\PP(N\ge b_t).
\end{eqnarray}
We shall show that $I_i=o(t^{-\alpha}L_1(t))$ for $i=1,3$ and $I_2\le cJ'_m\PP(X_1=t)$. 

The bound $I_3=o(t^{-\alpha}L_1(t))$ follows from
 (\ref{V-04-11}) and  (\ref{2019-01-07-1}) 
 with $\tau\ge 0$. 
 The bound $I_2\le cJ'_m\PP(X_1=t)$ follows from the 
relation
$t^{-1}\PP(X_1>t)(\alpha-1)\sim \PP(X_1=t)$,
the well known property of regularly varying sequences. 

Let us consider
$I_1$.
We denote
$a_t=t^2/(2(\alpha+1)\sigma^2\ln t)$~and~split
\begin{displaymath}
I_1
=
\left(\sum_{m< n< a_t}+\sum_{a_t\le n< b_t}\right)
n^{-1/2}e^{-t^2/(2n\sigma^2)}\PP(N=n)=:I_{1.1}+I_{1.2}.
\end{displaymath}
For $n\in (m,a_t)$ the inequality 
$e^{-t^2/(2n\sigma^2)}\le 
e^{-t^2/(2a_t\sigma^2)}=t^{-\alpha-1}$ implies
$I_{1.1}
=o(t^{-\alpha}L_1(t))$.
Let us show that  $I_{1.2}=o(t^{-\alpha}L_1(t))$.
We denote
$g(x)=x^{-1/2}e^{-t^2/(2x\sigma^2)}$ and 
$F(x)=-\PP(N\ge x)$ 
and apply the integration by parts formula (Sect. 2.9.24 in  \cite{Federer1969}) 
\begin{eqnarray}
I_{1.2}
=
\int_{a_t}^{b_t}g(x)dF(x)
=
F(b_t)g(b_t)-F(a_t)g(a_t)-\int_{a_t}^{b_t}g'(x)F(x)dx
=:
{\tilde J}_1+{\tilde J}_2+{\tilde J}_3.
\end{eqnarray}
It remains to show that 
${\tilde J}_i=o(t^{-\alpha}L(t))$, 
$i=1,2,3$. 
For $i=1,2$
the bounds are easy. 
We only prove the bound for ${\tilde J}_3$.
Note that 
$0<g'(x)
<
0.5(t/\sigma)^2x^{-5/2}e^{-t^2/(2x\sigma^2)}=:h(x)$, 
for $x\in(a_t,b_t)$. Invoking the latter inequality 
and  changing the variable of integration 
$z=t^2/(2x\sigma^2)$ we obtain
\begin{equation}\label{sausis4-3}
|{\tilde J}_3|
\le
\int_{a_t}^{b_t}h(x)\PP(N\ge x)dx
=
\sqrt{2}
\,
\frac{\sigma}{t}
\,
\int_{(1+\sigma^2)/(2\sigma^2)}^{(\alpha+1)\ln t}
\PP\bigl(N\ge t^2/(2\sigma^2z)\bigr)e^{-z}\sqrt{z}
\,
dz.
\end{equation} 
Given $\tau>0$ condition (\ref{2019-01-07-1}) implies 
that
$\PP(N\ge s^2\ln^{-\tau}s))
\le 
\phi(s)s^{1-\alpha}L_1(s)$ 
with some 
$\phi (s)\downarrow 0$ 
as $s\to+\infty$. 
Furthermore, from the inequality 
$u\ge t_u^2\ln^{-\tau}t_u$, with 
$t_u^2=4^{-\tau}u\ln^{\tau}u$,
which holds for sufficiently large $u$, we obtain
\begin{eqnarray}\label{2019-01-13+1}
\PP(N\ge u)
\le 
\PP(N\ge t_u^2\ln^{-\tau}t_u)
&
\le
&
\phi(t_u)t_u^{1-\alpha}L_1(t_u)
\\
\nonumber
&
\le
&
c\phi\bigl(\sqrt{4^{-\tau}u\ln^{\tau}u}\bigr) 
(u\ln^{\tau}u)^{(1-\alpha)/2}
L_1\bigl(\sqrt{u\ln^{\tau}u}\bigr).
\quad
\end{eqnarray}
Choosing $u=t^2/(2\sigma^2z)$ and 
invoking this inequality in (\ref{sausis4-3}) 
we obtain
\begin{equation}\label{sausis4-4}
|{\tilde J}_3|\le ct^{-\alpha}\phi(t^{1/2})
\int_{(1+\sigma^2)/(2\sigma^2)}^{(\alpha+1)\ln t}
e^{-z}z^{\alpha/2}
\ln^{(1-\alpha)\tau/2}(t^2/z)
L_1(tz^{-1/2}\ln^{\tau/2}t)dz.
\end{equation}
Here 
we estimated 
$\phi\bigl(\sqrt{4^{-\tau}u\ln^{\tau}u}\bigr)
\le \phi(\sqrt{t})$
using the monotonicity of $\phi$ and assuming that $t$ 
is sufficiently large.
 Finally, we show that the integral
 (\ref{sausis4-4}) is bounded from above by $cL_1(t)$.
Indeed,
for $z\in [(1+\sigma^2)/(2\sigma^2);(1+\alpha)\ln t]$ 
we have $\ln(t^2/z)\ge c'\ln t^2$. Furthermore, by 
general 
properties of slowly varying functions, see, e.g.,  
Th.1.1.4 
(iii) of \cite{Borovkov2008}, for any $\delta>0$ 
there exists $t_{\delta}>0$ such that 
$L_1(tz^{-1/2}\ln^{\tau/2}t)
\le (\ln t)^{\delta}L_1(t)$ for $t>t_{\delta}$.
Choosing $\delta<(\alpha-1)\tau/2$ we bound the
 integral 
 by  $cL_1(t)$. Hence 
 the right side of (\ref{sausis4-4})
 is $o(t^{-\alpha}L_1(t))$.
 
 Now assume that $\mu=0$ and (\ref{V-04-11+a}) holds.
 The only difference in the proof is that 
 for  $L_1(t)\sim c$ and 
 $\tau=0$ relation (\ref{2019-01-13+1}) reduces to
$\PP(N\ge u^2)\le \phi(u)u^{1-\alpha}$. Now 
 (\ref{sausis4-3}) implies
\begin{displaymath}
|{\tilde J}_3|
\le c
t^{-1}\phi(c't/\sqrt{\ln t})
\int_{(1+\sigma^2)/(2\sigma^2)}^{(\alpha+1)\ln t}
t^{1-\alpha}z^{\alpha/2}e^{-z}dz
=
o(t^{-\alpha}).
\end{displaymath}
The proof for $\mu=0$ is complete.

\medskip

 Assume that $\mu>0$. 
 We shall show that 
\begin{equation}\label{2019-01-11+1x}
I'_m
\le 
c'J'_m\PP(X_1=t)+o\bigl(\PP(X_1=t)\bigr).
\end{equation}
We put $t_*=t^{1/2}
w^{-1/2}\bigl(t/(2\mu)\bigr)$  in (\ref{VII-12-21}), then decompose 
$I'_m(t)$ using (\ref{VI-29-1}) and  derive (\ref{2019-01-11+1x})  from  the bounds shown below
\begin{eqnarray}\label{2019-01-11+3x}
&&
I'_{m.j}
\le 
c'\PP(X_1=t)w^{1/2}(t/(2\mu)),
\qquad
j=0,2,
\\
\label{2019-01-11+2x}
&&
I'_{m.1}
\le 
c'J'_m\PP(X_1=t),
\qquad
I'_{m.3}+I'_{m.4}
=o\bigl(\PP(X_1=t)\bigr).
\end{eqnarray}
Let us prove  (\ref{2019-01-11+3x}). The bound  for $I'_{m.0}$ is shown using the same argument as in (\ref{VII-30-3}) above. Next we consider $I'_{m.2}$.
For $n\in {\cal N}_2$ the inequalities $t/(2\mu)\le n\le t/\mu$ and (\ref{Borovkov8}) imply
\begin{equation}\label{2019-02-10-7}
\PP(S_n=t) =\PP(S_n^*=t_n^*) 
\le 
\frac{c_1}{\sqrt{t}}
e^{-c(t_n^*)^2/n}
+
c't(t_n^*)^{-\alpha}L_1(t_n^*)
=:a_n^{(1)}(t)+a_n^{(2)}(t).
\end{equation}
We show below that 
$A^{(k)}(t):=\sum_{n\in {\cal N}_2}a_n^{(k)}(t)
\le
cw^{1/2}(t/(2\mu))$, $k=1, 2$. These bounds combined with the inequalities that follow from (\ref{V-04-11}), (\ref{2017-06-04})
\begin{displaymath}
\PP(N=n)
\le
w(t/(2\mu))\PP(X_1=n)
\le 
c
w(t/(2\mu))\PP(X_1=t),
\qquad
n\in {\cal N}_2,
\end{displaymath}
imply $I'_{m.2}\le
c'\PP(X_1=t)w^{3/2}(t/(2\mu))$. Note that $w^{3/2}(\cdot)\le cw^{1/2}(\cdot)$, since $w(\cdot)\le c$.

Let us prove (\ref{2019-01-11+2x}).
For $n\in {\cal N}_1$ inequalities 
$t\ge t_n^*\ge t/2$, $t_n^*\ge n\mu$ and (\ref{Borovkov8}) imply 
\begin{equation}\label{2019-02-09+2}
\PP(S_n=t)=
\PP(S_n^*=t_n^*)
\le 
c' e^{-ct}
+
c'
n
t^{-\alpha}L_1(t).
\end{equation}
It follows that $\PP(S_n=t)\le c'nt^{-\alpha}L_1(t)$. The latter inequality implies
$I'_{m.1}\le cJ'_m\PP(X_1=t)$.
It remains to prove the second bound of (\ref{2019-01-11+2x}).
For $n\in {\cal N}_3\cup{\cal N}_4$ we have, cf. (\ref{2019-01-14+1}), 
\begin{eqnarray}\nonumber
\PP(S_n=t)
=
\PP({\hat S}_n={\hat t}_n)
&
\le 
&
\PP\bigl({\hat S}^{(1)}_n\ge{\hat t}_n/2,
\,
{\hat S}_n={\hat t}_n\bigr)
+
\PP\bigl({\hat S}^{(2)}_n\ge{\hat t}_n/2,
\,
{\hat S}_n={\hat t}_n\bigr)
\\
\label{2019-01-14+2}
&
\le
&
{\hat Q}^{(2)}_n\PP({\hat S}^{(1)}_n\ge{\hat t}_n/2)
+
{\hat Q}^{(1)}_n\PP({\hat S}^{(2)}_n\ge{\hat t}_n/2).
\end{eqnarray}
We estimate ${\hat Q}^{(k)}_n\le cn^{-1/2}$, $k=1,2$, by Lemma \ref{clem}.  
 Furthermore, we 
approximate the tail probabilities
$\PP({\hat S}_n^{(k)}\ge {\hat t}_n/2)$ by respective
Gaussian probabilities. The  non-uniform error
bound of \cite{Heyde}  implies
\begin{equation}\label{2019-01-14+4}
\PP({\hat S}_n^{(k)}\ge {\hat t}_n/2)
=
1-\Phi\bigl({\hat t}_n/(2\sigma_k))
+
c(1+({\hat t}_n/\sigma_k)^2)^{-1}\phi_n.
\end{equation}
Here 
$\sigma_k^2=\Var\bigl({\hat S}^{(k)}_n\bigr)$ 
and
$\phi_n
=
n^{-1/2}\EE|X_1|^3{\mathbb I}_{\{|X_1|\le \sqrt{n}\}}
+
\EE X_1^2{\mathbb I}_{\{|X_1|> \sqrt{n}\}}$.
Note that $\EE X_1^2<\infty$ implies
$\phi_n\to 0$ 
as $n\to+\infty$.
Using  relations $\sigma_k^2\asymp n$, $k=1,2$, and inequality 
$1-\Phi(x)\le e^{-x^2/2}$, $x>1$ we obtain from (\ref{2019-01-14+2}), (\ref{2019-01-14+4})
\begin{equation}\label{sausio3-5}
\PP(S_n=t)
\le 
c'n^{-1/2}
e^{-c{\hat t}^2_n/n}
+
c'n^{1/2}{\hat t}_n^{-2}\phi_n
=:b^{(1)}_n(t)+b^{(2)}_n(t).
\end{equation}
We show below that 
$B^{(k)}(t)
:=\sum_{n\in{\cal N}_3\cup{\cal N}_4}
b^{(k)}_n(t)
\le c'w^{1/2}(t/\mu)+c't^{-1/2}$, $k=1,2$.
These bounds together with (\ref{V-04-11}), (\ref{2017-06-04}) imply the second bound of 
(\ref{2019-01-11+2x}).

Now we evaluate the sums 
$A^{(k)}(t)$ and $B^{(k)}(t)$.
For $n\in {\cal N}_2$ inequalities 
$t_n^*\ge t-n\mu>0$ imply 
$(t_n^*)^2/n\ge \mu^3(n-t/\mu)^2/t$ and 
$(t_n^*)^3\ge \mu^3|n-t/\mu|^3$. We have
\begin{eqnarray}\label{2019-01-14+6}
&&
A^{(1)}(t)
\le 
\frac{c'}{\sqrt{t}}\sum_{n\in{\cal N}_2}
e^{-c(n-t/\mu)^2/t}
\le
\frac{c'}{\sqrt{t}}\int_{t_*}^{t/(2\mu)}e^{-cx^2/t}dx
\le 
c''\frac{\sqrt{t}}{t_*}
=
c''w^{1/2}\Bigl(\frac{t}{2\mu}\Bigr),
\qquad
\quad
\\
\label{2019-01-14+7}
&&
A^{(2)}(t)
\le c't\sum_{n\in{\cal N}_2}(t_n^*)^{-3}
\le
c't\int_{t_*}^{t/(2\mu)}x^{-3}dx
\le c'\frac{t}{t_*^2}
\le c'w\Bigl(\frac{t}{2\mu}\Bigr).
\end{eqnarray}
In (\ref{2019-01-14+6}) we used  inequality 
\begin{equation}\label{2019-01-14-8}
\int_{x\ge a}e^{-vx^2}dx
\le \int_{x\ge a}e^{-vx^2}\frac{vx}{va}dx
\le 
(2va)^{-1}
e^{-va^2}
\le (2va)^{-1},
\qquad
a,v>0.
\end{equation}
In (\ref{2019-01-14+7}) we applied  inequality 
$L_1(t_n^*)\le c_{\delta}(t_n^*)^{\delta}$ with 
$\delta=\alpha-3>0$, 
see Thm.1.1.4 in \cite{Borovkov2008}.

\noindent
Furthermore,
inequalities 
$b^{(1)}_n(t)\le c'n^{-1/2}e^{-cn}$,
$n\in {\cal N}_4$,
and 
$b^{(1)}_n(t)\le c't^{-1/2}e^{-c{\hat t}_n^2/t}$,
$n\in {\cal N}_3$, combined with
 (\ref{2019-01-14-8}) imply
\begin{eqnarray}
\label{2019-02-10-1}
&&
\sum_{n\in {\cal N}_4}b^{(1)}_n(t)
\le
c'
n^{-1/2}\sum_{n\in{\cal N}_4}e^{-cn}
\le c't^{-1/2}e^{-ct},
\\
\label{2019-02-10-2}
&&
\sum_{n\in {\cal N}_3}b^{(1)}_n(t)
\le 
\frac{c'}{\sqrt{t}}\int_{t^*}^{t/\mu}e^{-cx^2/t}dx
\le
c''\frac{\sqrt{t}}{t_*}=c''w^{1/2}\bigl(\frac{t}{2\mu}\Bigr).
\end{eqnarray}
Hence the bound for $B^{(1)}(t)$.  Finally, we bound $B^{(2)}(t)$ using the inequalities
\begin{eqnarray}
\label{2019-02-10-3}
\sum_{n\in {\cal N}_3}
\frac{\sqrt{n}}{{\hat t}_n^2}
\le
c'\sqrt{t}
\sum_{t_*\le j\le t/\mu}\frac{1}{j^2}
\le 
c''\frac{\sqrt{t}}{t_*},
\qquad
\sum_{n\in {\cal N}_4}
\frac{\sqrt{n}}{{\hat t}_n^2}
\le
c'\sum_{n\in{\cal N}_4}
\frac{1}{n^{3/2}}
\le 
\frac{c''}{\sqrt{t}}.
\end{eqnarray}
The proof for $\mu>0$ is complete.
$\qedsymbol$

{\it Proof of Theorem  \ref{R2019-16-1}}.
The proof is a straightforward extension of that of
Theorem 
\ref{T2+}, and we only indicate the changes.
In the proof we use  the property (A) of a slowly
 varying function 
$L$ that for any $\varepsilon>0$ there exists 
$t_{\varepsilon,L}>0$ such that 
$L(t)\le t^{\varepsilon}$ for 
$t>t_{\varepsilon,L}$.
Set $\delta=(2\alpha(\alpha-1)^{-1}+\tau)^{-1}$.

Let $1<\alpha\le 2$. We show that 
$I'_m(t)
\le 
c'\PP(X_1=t)\EE N^{1+\tau}{\mathbb I}_{\{N\ge m\}}$. 
Note that 
$\EE N^{1+\tau}{\mathbb I}_{\{N\ge m\}}=o(1)$ as 
$m\to+\infty$.
For this purpose we 
establish the bound 
\begin{equation}\label{2019-10-16+3}
\PP(S_n=t)\le cn^{1+\tau}t^{-\alpha}L_1(t).
\end{equation} 
 Lemma \ref{clem} and (\ref{Borovkov1}) imply
\begin{eqnarray}\nonumber
Q^{(k)}_n\le c'n^{-1/(\alpha-1)}/L_*(n),
\qquad
L^{(k)}_n(t,\delta)
\le
c' \bigl(n\PP(X_1\ge t\delta)\bigr)^{1/(2\delta)},
\qquad
k=1,2.
\end{eqnarray}
Furthermore, using property (A) we estimate
$Q^{(k)}_n\le c'n^{-1/(\alpha-1)}n^{\tau/2}$ and
\begin{displaymath}
L^{(k)}_n(t,\delta)
\le c''\bigl(nt^{1-\alpha}L_1(t)\bigr)^{1/(2\delta)}
\le 
cn^{\alpha/(\alpha-1)}
n^{\tau/2}
t^{-\alpha}L_1(t)
\end{displaymath}
so that the right side of (\ref{lema1+++}) is at most 
$c'''n^{1+\tau}t^{-\alpha}L_1(t)$.
Now Lemma \ref{lema1+} implies (\ref{2019-10-16+3}).

For $\alpha=2$ the same argument as in the case $1<\alpha\le 2$  above establishes 
(\ref{2019-10-16+3})
  for $n^{1+\tau/8}<t$. Indeed, we have that 
condition (\ref{Borovkov1+})  (that is required by (\ref{Borovkov1}) 
for $\alpha=2$) is implied by the inequality 
\begin{displaymath}
t^{1/(1+\tau/8)}
\PP(X_1>t\delta)
L_{\Delta}^{1+\beta}(t\delta)<\eta,
\end{displaymath}
where $\eta>0$ is
arbitrary, but independent of $t$ and $n$. The 
existence of such $\eta>0$ follows from the fact that,
by property (A), the left side 
is $o(1)$ as $t\to+\infty$.
For the remaining range $n^{1+\tau/8}\ge t$ we derive (\ref{2019-10-16+3}) from the concentration bound of Lemma \ref{clem} 
using property (A),
\begin{displaymath}
\PP(S_n=t)
\le 
 c\frac{n^{\tau/2}}{n^{1/(\alpha-1)}}
 \le
 c\frac{n^{\tau/2}}{n^{1/(\alpha-1)}}
 \Bigl(
 \frac{n^{2+\tau/4}}{t^2}
 \Bigr)^{\frac{8+2\tau}{8+\tau}}
 \le 
 c'n^{1+\tau}t^{-\alpha}L_1(t).
\end{displaymath}

\medskip

Let $2<\alpha<3$. For $\mu\le 0$ we have 
$t_n\ge t$.
 We derive (\ref{2019-10-16+3}) 
from inequalities (\ref{VII-11-2}), (\ref{Borovkov1}) 
using the same argument as in the case $1<\alpha<2$ 
above. We obtain 
\begin{equation}\label{2019-01-17+3}
\PP(S_n=t)=\PP({\tilde S}_n=t_n)
\le
cn^{1+\tau}t_n^{-\alpha}L_1(t_n)
\le
c'n^{1+\tau}t^{-\alpha}L_1(t).
\end{equation}

\noindent
Now assume that  $\mu>0$. 
We can assume that $\tau\le \beta$.
We put $t_*=t^{(\alpha-1)^{-1}+\beta/2}$
in (\ref{VII-12-21}), decompose 
$I'_m(t)$ using (\ref{VI-29-1}) and  estimate $I'_{m.j}$, $0\le j\le 4$. The concentration bound
$\PP(S_n=t)\le (n^{1/(\alpha-1)}L_*(n))^{-1}$ of Lemma \ref{clem}
together with $\PP(N=n)\le cn^{-\alpha-\beta}$ imply, cf. (\ref{VII-30-3}),
\begin{displaymath}
I'_{m.0}
\le 
c
t_*
\bigl(t^{1/(\alpha-1)}L_*(t)\bigr)^{-1}
t^{-\alpha-\beta}
\le 
ct^{-\alpha-\beta/2}(L_*(t))^{-1}
=
o\bigl(t^{-\alpha} L_1(t)\bigr).
\end{displaymath}
Furthermore, proceeding 
as in the proof of (\ref{2019-10-16+3})  above
we derive from (\ref{VII-11-2}), (\ref{Borovkov1}) the bounds  
\begin{eqnarray}\label{IX-08-1++}
&&
\PP(S_n=t)
=
\PP({\tilde S}_n=t_n)\le c'n^{1+\tau}t_n^{-\alpha}L_1(t_n),
\qquad
{\text{for}}
\qquad
n\in (m; t\mu^{-1}-t_*],
\\
\label{2019-01-10+1++}
&&
\PP(S_n=t)
=
\PP({\hat S}_n={\hat t}_n)
\le 
c'n^{1+\tau}{\hat t}_n^{-\alpha}L_1({\hat t}_n),
\qquad
{\text{for}}
\qquad
n\ge t\mu^{-1}+t_*.
\end{eqnarray}
These bounds imply 
$I'_{m.j}
\le  
c'\PP(X_1=t)\EE N^{1+\tau}{\mathbb I}_{\{N\ge m\}}$, 
$j=1,4$, see (\ref{IX-17-1}).
For $j=2,3$ inequalities (\ref{IX-08-1++}), (\ref{2019-01-10+1++}) 
imply
\begin{equation}\nonumber
\sum_{n\in{\cal N}_j}\PP(S_n=t)
\le
c't^{1+\tau}\sum_{n\in{\cal N}_j}|t-n\mu|^{-\alpha}L_1(|t-n\mu|)
\le
c'\frac{t^{1+\tau}}{t_*^{\alpha-1}}L_1(t^*).
\end{equation}
This bound combined with $\PP(N=n)\le ct^{-\alpha-\beta}$, $n\in {\cal N}_2\cup{\cal N}_3$,
implies $I'_{m.j}=o(t^{-\alpha}L_1(t))$.

\medskip

Let $\alpha=3$. We recall that for $\alpha=3$
 large deviation inequality (\ref{Borovkov1}) used in the proofs of (\ref{2019-01-17+3}) and (\ref{IX-08-1++}),
(\ref{2019-01-10+1++}) above only holds if additional condition (\ref{2019-01-10+2}) is satisfied. We write 
(\ref{2019-01-10+2}) in the form
\begin{equation}\label{2019-09-18}
 n
 \left(\frac{|\ln \Pi|}{y}\right)^2 
 L_{\diamond}\left(\frac{y}{|\ln\Pi|}\right)
\le \eta,
\qquad
{\text{where}}
\qquad
\Pi=\Pi(x)=nc_*x^{-2}L_1(x),
\end{equation}
and where $L_{\diamond}$ is a slowly varying function.

For $\mu<0$ we have $t_n\ge n|\mu|$. It is easy to show that for $x=t_n/2$ and $y=t_n\delta$  the left side of inequality
(\ref{2019-09-18})
is bounded uniformly in $n$ and 
$t$ (we assume that $n,t>A$, for some large  $A>0$). Now  (\ref{2019-01-17+3})
implies 
$I'_m(t)
\le 
c'\PP(X_1=t)\EE N^{1+\tau}{\mathbb I}_{\{N\ge m\}}$.

Assume that  $\mu=0$.  We apply (\ref{2019-09-16}) with
$a_t=t^{2-2\tau/(3+2\tau)}$. For $n<a_t$, $x=t/2$, $y=t\delta$  we  show that the left side of  inequality
(\ref{2019-09-18})
is bounded  using general properties of slowly varying functions and the fact that $n<t^{2-\varepsilon}$ for some given
$\varepsilon>0$. Then we apply (\ref{2019-01-17+3}) to 
$n<a_t$ and obtain the bound
$J_1
\le 
ct^{-\alpha}L_1(t)
\EE N^{1+\tau}{\mathbb I}_{\{m<N<a_t\}}$ 
cf. (\ref{2019-01-10+6}).
For $n\ge a_t$ we estimate 
$\PP(S_n=t)\le c'a_t^{-0.5}/L_*(a_t)$
using Lemma \ref{clem} and obtain the bound, 
cf. (\ref{2019-01-10+7}), 
\begin{displaymath}
J_2\le \frac{c'}{a_t^{0.5}L_*(a_t)}\PP(N>a_t)
\le\frac{c}{a_t^{1.5+\tau}L_*(a_t)}\EE N^{1+\tau}{\mathbb I}_{\{N>a_t\}}
=o(t^{-\alpha}L_1(t)).
\end{displaymath}
We obtain $I'_m(t)\le J_1+J_2\le ct^{-\alpha}L_1(t)
\EE N^{1+\tau}{\mathbb I}_{\{m<N<a_t\}} +o(t^{-\alpha}L_1(t))$.

Assume that  $\mu>0$. For $\alpha=3$ the proof is the same as that for $2<\alpha<3$ above. The only additional task is to check condition (\ref{2019-01-10+2}). Here we proceed similarly as in the proof of Theorem 2 (iv). 
%
$\qedsymbol$

	{\it Proof of Theorem \ref{R2019-16-2}}.
Note that  (\ref{2017-06-04}) implies $\EE N^{1+(\alpha-2)/2}<\infty$.
The proof of Theorem \ref{R2019-16-2}
is similar to that of Theorem
\ref{T2+} (iii).
Fix $A>0$. 
In (\ref{VII-12-21}) we put $t_*=Ab_t$, 
where
$b_t=t^{(\alpha-1)^{-1}}L_*(t)$ is defined in (\ref{2018-I-2}).
Then we decompose 
$I'_m(t)$ using (\ref{VI-29-1}) and  estimate $I'_{m.j}$, $0\le j\le 4$. We show that for some $c'>0$ (independent of $A$ and $t$) 
\begin{eqnarray}\label{2019-02-09-1}
&&
I'_{m.0}\le c'Aw(t/(2\mu))\PP(X_1=t),
\\
\label{2019-01-19+11}
&&
I'_{m.1}
\le c'\PP(X_1=t)\EE N^{1+0.5(\alpha-2)} {\mathbb I}_{\{N>m\}},
\\
\label{2019-01-19-2}
&&
I'_{m.j}
\le c'A^{1-\alpha}\PP(X_1=t),
\qquad
j=2,3,
\\
\label{2019-01-19+12}
&&
I'_{m.4}\le c' w(t/\mu)\PP(X_1=t).
\end{eqnarray}
These bounds together with (\ref{VI-29-1}) imply
$\limsup_t\frac{I'_m(t)}{\PP(X_1=t)}\le c'A^{1-\alpha}
+
c'\EE N^{1+0.5(\alpha-2)} {\mathbb I}_{\{N>m\}}$.
Letting $A,m\to+\infty$ we derive (\ref{V-12-3}) from
(\ref{VII-07+0}).

\noindent
Proof of (\ref{2019-02-09-1}). The bound 
$\PP(S_n=t)\le c'b^{-1}_n$ of Lemma \ref{clem}
together with (\ref{2017-06-04})  
imply (cf. (\ref{VII-30-3}))
\begin{displaymath}
I'_{m.0}
\le 
c'|{\cal N}_0|
\max_{n\in{\cal N}_0}
\Bigl\{
b^{-1}_n\PP(N=n)
\Bigr\}
\le c' Aw(t/2\mu)\PP(X_1=t).
\end{displaymath}
Proof of (\ref{2019-01-19-2}).  We only bound 
$I'_{m.2}$. The proof for $I'_{m.3}$ is much the same.
Using the bound ${\tilde Q}_n^{(k)}\le c'b^{-1}_n$,
$k=1,2$, of 
Lemma  \ref{clem} and invoking (\ref{Borovkov1}) we 
obtain from 
(\ref{VII-11-2}) for $\delta<0.5$ 
\begin{eqnarray}
&&
\PP({\tilde S}_n=t_n)
\le
c'n t_n^{-\alpha}L_1(t_n)
+
c_{\delta}b^{-1}_n\bigl(n\PP({\tilde X}_1>t_n\delta)\bigr)^{1/(2\delta)},
\qquad
n\in (m, t_{-}].
\qquad
\end{eqnarray}
Here $c_{\delta}>0$ only depends on $\delta$ and the distribution of $X_1$.  Put 
$\delta=(\alpha-1)/(2\alpha)$. We have
\begin{eqnarray}\label{2019-01-19-1}
\PP(S_n=t)
=
\PP({\tilde S}_n=t_n)
\le
c'n t_n^{-\alpha}L_1(t_n)
+
c'n t_n^{-\alpha}(L_1(t_n))^{\alpha/(\alpha-1)}(L_*(n))^{-1}.
\end{eqnarray}
For
$n\in (t/(2\mu), t_{-}]$ we estimate 
$n\le c't$ and 
$(L_*(n))^{-1}\le c'(L_*(t))^{-1}$ 
in  (\ref{2019-01-19-1}). Now Lemma \ref{lema3}
implies
\begin{equation}\label{2019-02-13-4}
\sum_{n\in{\cal N}_2}\PP(S_n=t)
\le 
c' A^{1-\alpha}.
\end{equation}
This bound combined with (\ref{2017-06-04}) shows 
$I'_{m.1}\le c'A^{1-\alpha}\PP(X_1=t)$.

\noindent
Proof of (\ref{2019-01-19+11}). 
We put
$\delta=(2\alpha(\alpha-1)^{-1}+2\tau)^{-1}$, 
where 
$\tau=(a-2)/4$. Now   (\ref{2019-01-19-1})  and 
$t/2\le t_n\le t$ imply for $n\in {\cal N}_1$
\begin{equation}\label{2019-01-19+13}
\PP(S_n=t)
\le 
c'nt^{-\alpha}L_1(t)
+
c'n^{1+\tau}(L_*(n))^{-1} t^{-\alpha-\tau(\alpha-1)}\bigl(L_1(t)\bigr)^{1+\tau+(\alpha-1)^{-1}}.
\end{equation}
For $\alpha>2$ inequality  $\tau>0$ implies 
$(L_*(n))^{-1}\le c'n^{\tau}$ and 
$(L_1(t))^{\tau+(\alpha-1)^{-1}}
\le 
c't^{\tau(\alpha-1)}$, by the 
general properties of slowly varying functions.
We obtain 
$\PP(S_n=t)\le cn^{1+2\tau}t^{-\alpha}L_1(t)$.
The latter inequality 
yields (\ref{2019-01-19+11}).

\noindent
Proof of (\ref{2019-01-19+12}). In view of
 (\ref{2017-06-04}) it suffices to show that 
$\sum_{n\in {\cal N}_4}\PP(S_n=t)\le c'$.
 Proceeding as in the proof of (\ref{2019-01-19-1}) we 
show  for $n\in {\cal N}_4$
\begin{displaymath}
\PP(S_n=t)
=
\PP({\hat S}_n={\hat t}_n)
\le
c'n {\hat t}_n^{-\alpha}L_1({\hat t}_n)
+
c'n {\hat t}_n^{-\alpha}
(L_1({\hat t}_n))^{\alpha/(\alpha-1)}
(L_*(n))^{-1}.
\end{displaymath}
Invoking the inequalities  
$n\mu\ge {\hat t}_n\ge n\mu/2$ 
we obtain
\begin{equation}\label{2019-02-13+1}
\PP(S_n=t)
\le 
cn^{1-\alpha}L_1(n)\bigl(1+(L_1(n))^{1/(\alpha-1)}/L_*(n)\bigr).
\end{equation}
We have $\sum_{n\in {\cal N}_4}\PP(S_n=t)
\le  c'$, since 
$\sum_n 
n^{1-\alpha}L_1(n)
\bigl(1+(L_1(n))^{1/(\alpha-1)}/L_*(n)\bigr)
<
\infty$.
\qquad\qquad
$\qedsymbol$

Now we prove Theorems \ref{T3+}, \ref{RR4} which establish 
 (\ref{VIII-27-3}), (\ref{2019-02-11-1}).
We begin with an outline of the proof. 
We 
split, see (\ref{b-17-1}),~(\ref{VI-29-1}),
\begin{equation}\label{XI-21-1}
\PP(S_N=t)
= 
I_m(t)+I'_m(t)
=
I_m(t)
+
I'_{m.0}+\dots+I'_{m.4},
\end{equation} 
where $I'_{m.j}$, $0\le j\le 4$ are defined by 
 (\ref{VII-12-21}) with 
$t_{\pm}=t\mu^{-1} \pm t_*$ and  $t_*:=u_t\mu^{-1}$.
The positive sequence $\{u_t\}$ will be
specified later. 
 For $\EE N<\infty$ we  will 
 choose $m=m_t\to+\infty$ as $t\to+\infty$ such that
\begin{equation}\label{2019-02-11++5}
I_m(t)=(\EE N)\PP(X_1=t)(1+o(1)).
\end{equation}
For $N$ satisfying $\EE N=\infty$ and (\ref{VIII-30-1}) 
with $\gamma>1$ we will 
choose $m=m_t\to+\infty$ as $t\to+\infty$ such that
\begin{equation}\label{2019-02-11+++2}
I_m(t)=o(\PP(N=t)).
\end{equation}
Furthermore, we show that as $t\to+\infty$
\begin{equation}\label{2019-02-11++1}
I'_{m.0}
= 
\mu^{-1}\PP\bigl(N=\lfloor t/\mu\rfloor\bigr)
(1+o(1))
\end{equation}
and that the remaining terms $I'_{m,i}$, $1\le i\le 4$,
on the right of (\ref{XI-21-1})
are negligibly small.

	{\it Proof of Theorem \ref{T3+}}.
In the proof we use the observation
that given a collection of sequences 
$\{a^{(k)}_t\}_{t\ge 1}$, $k=1,2,3,\dots$  
 such that $\forall k$ $\exists$ 
 $\lim_t a^{(k)}_t=:d_k$ and 
 $\sum_{k}|d_k|\PP(N=k)<\infty$,  one can find a 
 non-decreasing integer 
 sequence $m_t\to+\infty$ as 
$t\to+\infty$  such that 
\begin{equation}\label{IX-06-1}
\sum_{1\le k\le m_t}a^{(k)}_t\PP(N=k)
\to \sum_{k\ge 1}d_k\PP(N=k)
\qquad
{\text{as}}
\qquad
 t\to+\infty.
\end{equation}
\noindent
In (\ref{XI-21-1}) we choose 
$\{u_t\}$  satisfying 
(\ref{VIII-23-4})  
and  $u_t=o(t^{\varkappa})$ as $t\to+\infty$, 
where $\varkappa$ is from (\ref{IX-11-2}).

{\it Proof of} (\ref{2019-02-11-1}). The cases (i) and (ii) are treated simultaneously. Note that in both cases we require $\EE N<\infty$.
In the first step we establish 
(\ref{2019-02-11++5}) and (\ref{2019-02-11++1}).
 To show that there exists $\{m_t\}$ converging to $+\infty$ such that
(\ref{2019-02-11++5}) holds we apply
(\ref{IX-06-1}) to
$a^{(k)}_t:=\PP(S_k=t)/\PP(X_1=t)$ and use (\ref{V-04-10}) 
to verify 
the condition $\forall k$ 
$\lim_ta^{(k)}_t=n$.  
Next we derive  (\ref{2019-02-11++1}) from 
Lemma \ref{lema2++} 
  and (\ref{IX-11-2}).
In the second step we show that
\begin{eqnarray}
\label{2019-02-11++3}
&&
I'_{m.1}
=
c'\PP(X_1=t)
\EE
\bigl(N{\mathbb I}_{\{m_t\le N\le t/(2\mu)\}}\bigr),
\\
\label{2019-02-11++2}
&&
I'_{m.i}=o(\PP(N=t)),
\quad
i=2,3,
\\
\label{2019-02-11+++1}
&&
I'_{m.4}
=
o(\PP(X_1=t))
+
o(\PP(N=t)),
\quad
{\text{for}}
\quad 
\alpha>3. 
\qquad
\
\ 
\\
\label{2019-02-11++4}
&&
I'_{m.4}=o(\PP(X_1=t)),
\quad
{\text{for}}
\quad 
\alpha\le 3,
\end{eqnarray}
Note that (\ref{2019-02-11++3}) combined with 
$\EE N<\infty$ implies $I'_{m.1}=o(\PP(X_1=t))$.
Invoking (\ref{2019-02-11++5}), (\ref{2019-02-11++1}), (\ref{2019-02-11++3}-\ref{2019-02-11++4})
 in 
(\ref{XI-21-1}) we obtain (\ref{2019-02-11-1}).
It remains to prove 
(\ref{2019-02-11++3}-\ref{2019-02-11++4}).

\noindent
Proof of (\ref{2019-02-11++3}).
  For $\alpha\le 3$, respectively $\alpha>3$,
   we obtain from 
(\ref{IX-08-1}), (\ref{IX-17-1}), respectively
(\ref{2019-02-09+2}), that 
$\PP(S_n=t)\le c'n\PP(X_1=t)$, $n\in {\cal N}_1$.
This bound implies (\ref{2019-02-11++3}).

\noindent
Proof of (\ref{2019-02-11++2}) 
for $\alpha\le 3$. 
From (\ref{IX-18-1++}) 
we obtain for $i=2,3$ that
\begin{displaymath}
\sum_{n\in{\cal N}_i}\PP(S_n=t)
\le 
c'
t\, t_*^{1-\alpha}
\le
c''t\, u_t^{1-\alpha}.
\end{displaymath}
This bound together with (\ref{IX-11-1}) imply
$I'_{m.i}\le c'tu_t^{1-\alpha}\PP(N=t)=
o\bigl(\PP(N=t)\bigr)$.
In the last step we used 
$t=o(u_t^{\alpha-1})$, see  the 
first relation of (\ref{VIII-23-4}).

\noindent
Proof of (\ref{2019-02-11++2}) for $\alpha>3$ and (\ref{2019-02-11+++1}).
From (\ref{sausio3-5}), (\ref{2019-02-10-1}),
(\ref{2019-02-10-2}), (\ref{2019-02-10-3}) we obtain
\begin{eqnarray}
\label{2019-02-10-4}
&&
I'_{m.3}
\le 
\max_{n\in{\cal N}_3}\PP(N=n)
\sum_{n\in{\cal N}_3}\PP(S_n=t)
\le 
\max_{n\in{\cal N}_3}\PP(N=n)
ct^{1/2}t_*^{-1},
\\
\label{2019-02-10-5}
&&
I'_{m.4}
\le
\sum_{n\in{\cal N}_4}b^{(1)}(t)+
\sum_{n\in{\cal N}_4}b^{(2)}(t)\PP(N=n)
\le 
c'e^{-ct}
+
\sum_{n\in{\cal N}_4}n^{-3/2}\PP(N=n).
\qquad
\end{eqnarray}
The right side of (\ref{2019-02-10-4}) is 
$o(\PP(N=t))$
since
$\max_{n\in{\cal N}_3}\PP(N=n)\le c\PP(N=t)$ and 
$\sqrt{t}t_*^{-1}=o(1)$ by (\ref{IX-11-1}) and 
(\ref{VIII-23-4}) respectively. The 
right side of (\ref{2019-02-10-5}) is 
$o(\PP(X_1=t))+o(\PP(N=t))$
by (\ref{2019-02-10-6}) and 
 because $e^{-ct}=o(\PP(X=t))$.
Similarly,
from (\ref{2019-02-10-7}), (\ref{2019-01-14+6}),
(\ref{2019-01-14+7}) we obtain
\begin{displaymath}
I'_{m.2}
\le 
\max_{n\in{\cal N}_2}\PP(N=n)
\sum_{n\in{\cal N}_2}\PP(S_n=t)
\le 
\max_{n\in{\cal N}_3}\PP(N=n)
c(t^{1/2}t_*^{-1}+tt_*^{-2})
=
o(\PP(N=t)).
\end{displaymath}

\noindent
Proof of (\ref{2019-02-11++4}). 
From 
(\ref{2019-01-10+1}), (\ref{IX-17-1})  we have
$\PP(S_n=t)\le c'nt^{-\alpha}$, $n\in {\cal N}_4$.
This bound implies 
$I'_{m.4}\le c'\PP(X_1=t)
\EE (N{\mathbb I}_{\{N\ge 2t/\mu\}})$. But 
$\EE (N{\mathbb I}_{\{N\ge 2t/\mu\}})=o(1)$ since 
$\EE N<\infty$.

{\it Proof of} (\ref{VIII-27-3}).
The cases (i) and (ii) are treated simultaneously.
For $\EE N<\infty$ we derive (\ref{VIII-27-3})
 from  (\ref{VIII-27-2}), (\ref{2019-02-11-1}).
It remains to consider the case where (\ref{VIII-30-1}) 
holds and $\EE N=\infty$.
To show that there exists $\{m_t\}$ converging to $+\infty$ such that (\ref{2019-02-11+++2}) holds
we apply  
(\ref{IX-06-1}) to
$a^{(k)}_t:=\PP(S_k=t)/\PP(N=t)$ and use 
(\ref{V-04-10}) and
(\ref{VIII-27-2}) to verify 
the condition $\forall k$ 
$\lim_ta^{(k)}_t=0$. We remark that 
(\ref{2019-02-11++1}), (\ref{2019-02-11++3}),
(\ref{2019-02-11++2}), 
(\ref{2019-02-11+++1})
remain true as their proofs above have not used the 
condition $\EE N<\infty$.
Next we show that 
\begin{equation}\label{2019-02-11+++3}
I'_{m.1}=o(\PP(N=t))
\quad
{\text{for}}
\quad
\alpha>2
\quad
{\text{and}}
\quad
I'_{m.4}=o(\PP(N=t))
\quad
{\text{for}}
\quad
\alpha\le 3.
\end{equation}
Note that 
(\ref{2019-02-11+++2},\ref{2019-02-11++1},\ref{2019-02-11++2},\ref{2019-02-11+++1},\ref{2019-02-11+++3})
combined with
(\ref{XI-21-1}) yield
(\ref{VIII-27-3}).
Hence it remains to prove (\ref{2019-02-11+++3}).
Note that
(\ref{VIII-30-1})  and $\EE N=\infty$
imply 
$\EE N {\mathbb I}_{\{m_t\le N\le t/(2\mu)\}}$
$\le $
$\EE N {\mathbb I}_{\{N\le t/(2\mu)\}}$
$\le c' 
t^{2-\gamma}{\bar L}(t)$, where ${\bar L}$ is a slowly varying 
function.
 Invoking this inequality in  (\ref{2019-02-11++3}) we 
 obtain
the first bound of (\ref{2019-02-11+++3}).
Let us show that (\ref{VIII-30-1}) 
implies the second bound of 
(\ref{2019-02-11+++3}).
Invoking the inequalities
$\PP(S_n=t)\le cn^{-1/(\alpha-1)}$ for $n>t^{\alpha-1}$ and $\PP(S_n=t)\le c'nt^{-\alpha}$ for $n\le t^{\alpha-1}$
 (the first one follows from Lemma \ref{clem}, the second one follows from (\ref{2019-01-10+1}), (\ref{IX-17-1})) we obtain
\begin{eqnarray}\nonumber
I'_{m.4}
&
\le
&
c't^{-\alpha}\sum_{2t/\mu\le n\le t^{\alpha-1}}n^{1-\gamma}L_2(n)
+
c'\sum_{n>t^{\alpha-1}}n^{-\gamma -(\alpha-1)^{-1}}L_2(n)
\\
\nonumber
&
\le
&
c't^{-\alpha}\bigl(t^{\alpha-1}\bigr)^{2-\gamma}L(t^{\alpha-1})
+
c'\bigl( t^{\alpha-1}\bigr)^{1-\gamma-(\alpha-1)^{-1}} L_2(t^{\alpha-1})
\\
\nonumber
&
=
&
c't^{-\gamma-(\alpha-2)(\gamma-1)}\Bigl( L(t^{\alpha-1})+L_2(t^{\alpha-1})\bigr)
\\
\nonumber
&
=
&
o\bigl(\PP(N=t)\bigr).
\end{eqnarray}
Here $L$ is a slowly varying function ($L=L_2$, for $\gamma\not=2$).
In the last step we used $(\alpha-2)(\gamma-1)>0$. 
$\qedsymbol$

{\it Proof of Theorem \ref{RR4}}.
The proof is similar to that of
Theorem \ref{T3+} (i).
We only prove (\ref{2019-02-11-1}). The proof of 
(\ref{VIII-27-3}) is almost the same.

{\it Proof of} (\ref{2019-02-11-1}).
We decompose $\PP(S_N=t)$ using
(\ref{XI-21-1}) and 
(\ref{VII-12-21}), where we put 
$t_{\pm}=t\mu^{-1} \pm t_*$, $t_*:=Ab_t$ and where 
$A>0$ is a large number and $b_t$ is from (\ref{2018-I-2}). 
We choose $m=m_t\to+\infty$ such that
(\ref{2019-02-11++5}) holds.
Furthermore, we show below  that for some $c'>0$ 
(independent of $A$ and~$t$) 
\begin{eqnarray}\label{2019-02-13-1}
&&
I'_{m.1}
\le 
c'\PP(X_1=t)
\EE N {\mathbb I}_{\{m\le N\le t/(2\mu)\}},
\\
\label{2019-02-13-2}
&&
I'_{m.j}
\le c'A^{1-\alpha}\PP(N=t),
\qquad
j=2,3,
\\
\label{2019-02-13-3}
&&
I'_{m.4}=o(\PP(N=t)).
\end{eqnarray}
Finally, we prove that 
\begin{equation}\label{XI-22-2}
\limsup_{t\to+\infty}\Bigr|
\frac{I'_{m.0}}
{\PP(N=\lfloor t/\mu\rfloor)}
-
\mu^{-1}
\Bigr|
\to 0
\qquad
{\text{as}}
\qquad A\to+\infty.
\end{equation}
Invoking (\ref{2019-02-11++5}), (\ref{2019-02-13-1}),
(\ref{2019-02-13-2}), (\ref{2019-02-13-3}),
(\ref{XI-22-2})
in (\ref{XI-21-1})
and using our assumption $\EE N<\infty$  we obtain 
(\ref{2019-02-11-1}), by  
letting $A\to+\infty$.

\noindent
Proof of (\ref{2019-02-13-1}).
It suffices to show that
$\PP(S_n=t)\le c'nt^{-\alpha}L_1(t)$.
We choose  $\delta<(\alpha-1)/(2\alpha)$ and apply 
(\ref{VII-11-2}) to $\PP(S_n=t)=\PP({\tilde S}_n=t_n)$. 
Invoking in (\ref{VII-11-2}) the inequalities
\begin{displaymath}
\max_{\delta t_n\le j\le t_n}
\PP({\tilde X}_1=j)
\le
c't_n^{-\alpha}L_1(t_n),
\quad
\
{\tilde Q}_n^{(k)}\le b_n^{-1},
\quad
\
c'{\tilde L}_n^{(k)}(t_n,\delta)
\le
c'
\bigl(nt_n^{1-\alpha}L_1(t_n)\bigl)^{1/(2\delta)} 
\end{displaymath}
(the last one follows (\ref{Borovkov1})) and using 
$t/2\le t_n\le t$, for $n\in {\cal N}_1$, we obtain
\begin{eqnarray}
\nonumber
\PP({\tilde S}_n=t_n)
&
\le
& 
c'n t_n^{-\alpha}L_1(t_n)
+
c'n^{(2\delta)^{-1}-(\alpha-1)^{-1}}
(L_*(n))^{-1}
t_n^{(1-\alpha)/(2\delta)}(L_1(t_n))^{(2\delta)^{-1}}
\\
\label{XII-19-1}
&
\le
& 
c'n t^{-\alpha}L_1(t)
+
c'n^{(2\delta)^{-1}-(\alpha-1)^{-1}}
(L_*(n))^{-1}
t^{(1-\alpha)/(2\delta)}(L_1(t))^{(2\delta)^{-1}}
\\
\nonumber
&
\le 
&
c'nt^{-\alpha}L_1(t).
\end{eqnarray}
To prove the last inequality
 we write
the second summand on the right of (\ref{XII-19-1}) in 
the form
 \begin{displaymath}
 c'n t^{-\alpha}L_1(t)R_n(t),
 \qquad
 R_n(t):=
 \bigl(
 n^{1/(\alpha-1)}t^{-1}
 \bigr)^{\alpha\tau}
 (L_*(n))^{-1}
 (L_1(t))^{\alpha(1+\tau)(\alpha-1)^{-1}-1}
 \end{displaymath}
 and observe that $R_n(t)$ is bounded uniformly 
 in $n\in {\cal N}_1$.
Here $\tau>0$ is defined by the 
equation $1/(2\delta)=(\alpha/(\alpha-1))(1+\tau)$.
 Indeed,
the inequality $n\le t/(2\mu)$ (which holds for
$n\in{\cal N}_1$) implies  $n^{1/(\alpha-1)}t^{-1}\le c't^{-\varepsilon'}$ with $\varepsilon'=1-(\alpha-1)^{-1}>0$. In addition,
by the properties of slowly varying functions, we have 
  $|L^{-1}_*(n)|=o( n^{\varepsilon})$ and  
 $|L_1(t)|= o( t^{\varepsilon})$ for any $\varepsilon>0$ as $n,t\to+\infty$. Hence, $R_n(t)\le c'$ uniformly in $n\in{\cal N}_1$.
 
\noindent
Proof of (\ref{2019-02-13-2}). This bound follows 
from  (\ref{IX-11-1})  and (\ref{2019-02-13-4}).

\noindent
Proof of (\ref{2019-02-13-3}). Using general properties 
of slowly varying functions we obtain from 
(\ref{2019-02-13+1}) that 
$\PP(S_n=t)\le c'n^{-1-\beta}$, $n\in {\cal N}_4$. 
Now (\ref{2019-02-13-3}) follows from (\ref{IX-11-1}) 
and (\ref{2019-02-10-6}).

\noindent
Proof of (\ref{XI-22-2}).
Proceeding as in the proof of Lemma \ref{lema2++} we show that
 for any (small) $\delta>0$ and (large) $A_0>0$ one can 
find  $A>A_0$ and large $t_0$ such that 
\begin{equation}\label{XI-22-1}
 \forall
 \,
  t>t_0
\qquad
 \Bigl|\sum_{n\in {\cal N}_0}\PP(S_n=t)-\mu^{-1}\Bigr|<\delta.
\end{equation}
We obtain (\ref{XI-22-2}) from  (\ref{XI-22-1}) and the relation that follows from (\ref{IX-11-2})  
\begin{equation}\nonumber
\qquad\qquad\qquad
\max_{n\in {\cal N}_0}
\Bigl|\frac{\PP(N=n)}{\PP(N=\lfloor t/\mu\rfloor)}-1\Bigr|
\
\to 
0
\qquad
{\text{ as}}
\qquad
 t\to+\infty.
 \qquad
 \qquad\qquad
 \qquad\qquad\qquad
 \qedsymbol
 \end{equation}



	{\it Proof of Theorem \ref{T1+}}.
(i) follows from
Theorem \ref{T2+} (v) and Theorem \ref{R2019-16-1}
for $\alpha>3$ and $\alpha\le 3$ respectively.

\noindent
Proof of (ii) and (iii). Given $1< m_t <t/(2\mu)$ we split
\begin{equation}\label{2019-09-08}
\PP(S_N=t)
=
\EE
\bigl(\PP(S_N=t|N)
\bigl(
{\mathbb I}_{\{N\le m_t\}}
+
{\mathbb I}_{\{m_t<N<\frac{t}{2\mu}\}}
+
{\mathbb I}_{\{\frac{t}{2\mu}\le N\}}
\bigr)
\bigr)
=:{\tilde I}_1+{\tilde I}_2+{\tilde I}_3.
\end{equation}
We first 
evaluate ${\tilde I}_3$.
We choose a nondecreasing sequence $u_t\to+\infty$ which satisfies (\ref{VIII-23-4}) and put
$t_{\pm}=(t\pm u_t)\mu^{-1}$. It follows from 
Lemma \ref{lema2++} and the renewal theorem \cite{ErdosFellerPollard1949} that 
\begin{displaymath}
\sum_{t_{-}\le n\le t_{+}}\PP(S_n=t)\sim \mu^{-1} 
\qquad
{\text{and}}
\qquad
\sum_{n\ge 1}\PP(S_n=t)\sim \mu^{-1}
\qquad {\text{as}}
\qquad 
t\to+\infty.
\end{displaymath}
In particular, we have
$\sum_{n\ge 1, n\notin[t_{-}, t_{+}]}\PP(S_n=t)=o(1)$.
In view of (\ref{VIII-30-1}) these relations imply
\begin{equation}\label{2019-09-09}
{\tilde I}_3=
\EE 
\bigl(
\PP(S_N=t|N){\mathbb I}_{\{t_{-}\le N\le t_{+}\}}
\bigr)
+
o(\PP(N=t))
\sim
\mu^{-1}\PP(N=\lfloor t/\mu\rfloor).
\end{equation}
Next, we estimate 
${\tilde I}_2$. By Theorem 1 of \cite{Doney1989}, for any sequence $\{m_t\}$ such that $m_t\uparrow +\infty$ as $t\to+\infty$ there is a constant $c>0$ such that
$\PP(S_n=t)\le c n\PP(X_1=\lfloor t - n\mu\rfloor)$  
for  $m_t\le n\le t/(2\mu)$. In view of (\ref{V-04-11}) we obtain for some constant $c'>0$ (depending on $\{m_t\}$) that
\begin{displaymath}
{\tilde I}_2
\le 
c'
\PP(X_1=t)\xi_t,
\qquad \xi_t:=
\EE
\bigl( 
N {\mathbb I}_{\{m_t\le N\le \frac{t}{2\mu}\}}\bigr).
\end{displaymath}

Let us prove (iii).
We choose $m_t\uparrow+\infty$ as $t\to+\infty$ such that 
$
{\tilde I}_1\sim (\EE N)\PP(X_1=t)$, see (\ref{2019-02-11++5}).
Furthermore, 
for $\EE N<\infty$ we have 
$\xi_t\le \EE (N{\mathbb I}_{\{N\ge m_t\}})=o(1)$. Hence ${\tilde I}_2=o(\PP(X_1=t))$.
Collecting these relations and (\ref{2019-09-09}) in 
(\ref{2019-09-08}) we obtain (iii).

Let us prove (ii). We only consider the case 
$\EE N=\infty$. 
We choose $m_t\uparrow+\infty$ such that
${\tilde I}_1=o(\PP(N=t))$, see (\ref{2019-02-11+++2}).  
Furthermore, for $\EE N=\infty$
we have, by  (\ref{VIII-30-1}),
$\xi_t\le \EE
\bigl( 
N {\mathbb I}_{\{N\le \frac{t}{2\mu}\}}\bigr)
= t^{2-\gamma}{\bar L}_2(t)\to+\infty$ as 
$t\to+\infty$, where ${\bar L}_2$ is slowly varying at infinity. 
Hence 
${\tilde I}_2
\le 
ct^{2-\alpha-\gamma}L_1(t){\bar L}_2(t)
=
o(\PP(N=t))$. 
Invoking these bounds and (\ref{2019-09-09}) in 
(\ref{2019-09-08}) we obtain (ii).

Proof of (iv). For $n\to+\infty$ the 
standardized sums $n^{-1/(\alpha-1)}(X_1+\dots+X_n)$ converge 
in distribution to an $\alpha-1$ stable random variable, which 
we denote by $Z_a$. Here the subscript $a$ refers to the 
constant $a$ in (\ref{V-04-11+a}). Note that $Z_a$ and $a^{1/\alpha}Z_1$ have the same  distributions. Therefore, it suffices to show that
\begin{equation}\label{2017-XII-29-1}
\PP(S_N=t)
\sim
h(t)
\EE Z_a^{(\alpha-1)(\gamma-1)},
\qquad
h(t)
:=
t^{-1-(\alpha-1)(\gamma-1)}L_2(t^{\alpha-1})(\alpha-1).
\end{equation}
Given $A>0$ denote 
$J_A=\EE \bigl( Z_a^{(\alpha-1)(\gamma-1)}{\mathbb I}_{\{A^{-1}\le Z_a^{\alpha-1}\le A\}}\bigr)$.
We prove below that 
\begin{equation}\label{2017-XII-29-2}
J_A
\le 
\liminf_t \frac{\PP(S_N=t)}{h(t)}
\le 
\limsup_t
\frac{\PP(S_N=t)}{h(t)}
\le J_A+c'\bigl(A^{1-(\alpha-1)^{-1}-\gamma}+A^{\gamma-2}\bigr).
\end{equation}
Then 
(\ref{2017-XII-29-1}) follows from
(\ref{2017-XII-29-2})
by letting $A\to+\infty$.
Let us prove (\ref{2017-XII-29-2}).
We split
\begin{eqnarray}\nonumber
&&
\PP(S_N=t)
=
\sum_{n\ge 1}\PP(S_n=t)\PP(N=n)
= I^*_1+I^*_2+I^*_3,
\qquad
I^*_j=\sum_{n\in {\cal N}_j}\PP(S_n=t)\PP(N=n)
\\
\nonumber
&&
{\cal N}_1
=
\{n\le A^{-1}t^{\alpha-1}\},
\qquad
{\cal N}_2
=
\{A^{-1}t^{\alpha-1}
<n
< 
At^{\alpha-1}\},
\qquad
{\cal N}_3
=
\{n\ge At^{\alpha-1}\}.
\end{eqnarray}
We first show that as $t\to+\infty$
\begin{equation}\label{2019-09-19++}
I^*_2=h(t)(J_A+o(1)).
\end{equation}
 From the local limit theorem bound (\ref{VIII-23-3}) we obtain 
\begin{eqnarray}\nonumber
&&
I^*_2
=
\sum_{n\in{\cal N}_2}
(g(t_n)+\delta_n)b_n^{-1}
\PP(N=n)
=
I_2^{**}+R,
\\
\nonumber
&&
I_2^{**}:=\sum_{n\in{\cal N}_2}
g(t_n)b_n^{-1}
\PP(N=n),
\qquad
R:=\sum_{n\in{\cal N}_2}
\delta_nb_n^{-1}
\PP(N=n).
\end{eqnarray}
Here $t_n=tb_n^{-1}$ and $b_n=n^{1/(\alpha-1)}$.
  $g(\cdot)$ stands for  the density of $Z_a$ and $\delta_n$ denotes the remainder term. We have $|\delta_n|\le \tau_n$ and $\max_{k\in{\cal N}_2}|\delta_k|\to 0$ as $t\to+\infty$. Recall that $\tau_n$ is defined in (\ref{VIII-23-3}).
   Using the relation
\begin{equation}\nonumber
t_n-t_{n+1}=t_n(\alpha-1)^{-1}n^{-1}\bigl(1+O(n^{-1})\bigr),
\end{equation}
we write $I^{**}_2$ in the form $I^{**}_2=h(t)S$,
where
\begin{equation}\nonumber
S=\sum_{n\in{\cal N}_2}
g(t_n)t_n^{(\alpha-1)(\gamma-1)}(t_n-t_{n+1})
\frac{1}{1+O(n^{-1})}\frac{L_2(n)}{L_2(t^{\alpha-1})}
\end{equation}
converges to $J_A$ as $t\to+\infty$. Here we used the fact that $\frac{1}{1+O(n^{-1})}\frac{L_2(n)}{L_2(t^{\alpha-1})}\to 1$ uniformly in $n\in{\cal N}_2$. 
Now (\ref{2019-09-19++}) follows from the simple bound
\begin{displaymath}
R
\le
\bigl(\max_{k\in{\cal N}_2}|\delta_k|\bigr)\sum_{n\in{\cal N}_2} b_n^{-1}\PP(N=n)
=o(1)\sum_{n\in{\cal N}_2} b_n^{-1}\PP(N=n)
=o(h(t))
\quad
{\text{as}}
 \quad
 t\to+\infty.
\end{displaymath}
We secondly  estimate $I^*_j$, $j=1,3$.  Using the  
the local limit theorem bound 
$\PP(S_n=t)\le c'n^{-1/(\alpha-1)}$ for $n\in {\cal N}_3$ and 
the bound $\PP(S_n=t)\le c'n t^{-\alpha}$ for $n\in{\cal N}_1$, see (\ref{VII-1}), we obtain as $t\to+\infty$
\begin{eqnarray}
\nonumber
\qquad
I^*_3
&
\le 
&
c'\sum_{n\in{\cal N}_3}
n^{-1/(\alpha-1)}\PP(N=n)
\sim
c't^{-1-(\alpha-1)(\gamma-1)}L_2(t^{\alpha-1})
A^{1-(\alpha-1)^{-1}-\gamma},
\\
\nonumber
I^*_1
&
\le 
&
c't^{-\alpha}\sum_{n\in{\cal N}_1}n\PP(N=n)
\sim
c't^{-1-(\alpha-1)(\gamma-1)}L_2(t^{\alpha-1})A^{\gamma-2}.
\end{eqnarray}
%
These relations together with (\ref{2019-09-19++}) imply (\ref{2017-XII-29-2}).
\quad
$\qedsymbol$

	{\it Proof of Proposition \ref{X-12-1}}.
The proof of statement (i) is much the same as that of  
Theorem \ref{T3+} (i), but now 
in (\ref{XI-21-1}) we choose 
$u_t=\bigl(t\ln t\bigr)^{1/2}\ln\ln t$.
Note that $\{u_t\}$ satisfies 
(\ref{VIII-23-4}) (to check the third condition of 
(\ref{VIII-23-4}) we use Remark \ref{X-11-1} and 
(\ref{alpha=3})).

To show 
(ii) we proceed similarly as in the proof of Theorem \ref{T2+} (iv) for $\mu>0$. The only difference is that now instead of (\ref{2019-01-11+3})
we use the bound that follows from Lemma \ref{lema2++}
\begin{displaymath}
\qquad
\qquad
I'_{m.0}\le \max_{n\in{\cal N}_0}\PP(N=n)
\sum_{n\in{\cal N}_0}\PP(S_n=t)
\le 
w(t/(2\mu))t^{-3}(\mu^{-1}+o(1))=o(t^{-3}).
\qquad
\qquad\qquad
\qedsymbol
\end{displaymath}

	{\it Proof of relation (\ref{VIII-6-1})}.
For deterministic $n$ relation
$\PP(M_n=t)\sim n\PP(X_1=t)$ follows from the inequalities
\begin{equation}\label{V-11-1}
np^*-\binom{n}{2}p^{**}
\le 
\PP(M_n=t)\le
np^*,
\end{equation}
where 
\begin{eqnarray}\nonumber
p^*
&
=
&
\PP(X_n=t, M_{n-1}\le t)=\PP(X_n=t)\PP(M_{n-1}\le t)\sim \PP(X_1=t),
\\
\nonumber
p^{**}
&
=
&
\PP(X_1=t, X_2=t)=\PP(X_1=t)\PP(X_2=t) =o\bigl(\PP(X_1=t)\bigr).  
\end{eqnarray}

Let us prove (\ref{VIII-6-1}). To this aim we show that for any $\varepsilon>0$
\begin{equation}\label{V-11-2}
 (1-\varepsilon)(\EE N)
\le 
\liminf_{t\to+\infty} \frac{\PP(M_N=t)}{\PP(X_1=t)}
\le
\limsup_{t\to+\infty} \frac{\PP(M_N=t)}{\PP(X_1=t)}
\le
\EE N.
\end{equation}
To show the very left inequality we choose large positive integer $m$ such that
$\EE N{\mathbb I}_{\{N\le m\}}>(1-\varepsilon)\EE N$ and 
use the left inequality of  (\ref{V-11-1}). We obtain
\begin{eqnarray}
 \nonumber
\PP(M_N=t)
&
\ge 
&
\EE \Bigl(\PP(M_N=t|N){\mathbb I}_{\{N\le m\}}\bigr)
\\
\nonumber
&
\ge
& 
(1+o(1))\PP(X_1=t)\EE N{\mathbb I}_{\{N\le m\}}
+o\bigl(\PP(X_1=t)\bigr).
\end{eqnarray}
The very right inequality of (\ref{V-11-2}) follows by  
Lebesgue's dominated convergence theorem from the right 
inequality of (\ref{V-11-1})
\begin{displaymath}
\qquad
\PP(M_N=t)
=
\EE \bigl(\PP(M_N=t|N)\bigr)
\le
\EE \bigl(N\PP(X_1=t)\bigr)
=
(\EE N)\PP(X_1=t).
\qquad
\qquad
\qedsymbol
\end{displaymath}
%


	{\it Proof of Theorem \ref{T4+}}.
Before the proof we introduce some notation.  Denote  $a_i=\EE X_1^i$, $b_i=\EE Y_1^i$.
 For $i=0,1$ we denote by $\Lambda_i^{(r)}$
a mixed Poisson random variable with the distribution
\begin{displaymath}
\PP(\Lambda_i^{(r)}=s)
=
\bigl(
\EE \lambda^r_i
\bigr)^{-1}
\EE\Bigl( e^{-\lambda_i}\lambda_i^{s+r}/s!\Bigr),
\qquad
s=0,1,2,\dots.
\end{displaymath}
Here  $\lambda_0=Y_1\beta^{1/2}a_1$ and $\lambda_1=X_1\beta^{-1/2}b_1$ . 
 Let $\tau_1,\tau_2,\dots$ be iid copies of $\Lambda_1^{(1)}$. 
Assuming that 
$\{\tau_i$, $i\ge 1\}$ are independent of 
$\Lambda_0^{(r)}$, define randomly stopped sums 
\begin{displaymath}
d_*^{(r)}=\sum_{j=1}^{\Lambda_0^{(r)}}\tau_j,
\qquad
r=0,1,2.
\end{displaymath}
Finally, we denote 
\begin{equation}
\label{X-28-2+}
C_*(k)
=
\left(
1+\sqrt{\beta} \, \frac{a_2^2b_2}{a_3b_1}\frac{p_1(k)}{p_2(k)}
\right)^{-1}.
\end{equation}
Here 
\begin{equation}\label{X-30-1}
p_1(k)
=
\PP\bigl(d_*^{(2)}+\Lambda_1^{(2)}+{\bar \Lambda}_1^{(2)}=k-2\bigr),
\qquad
p_2(k)
=
\PP\bigl(d_*^{(1)}+\Lambda_1^{(3)}=k-2\bigr).
\qquad
\end{equation}
The random variables $d_*^{(1)}, \Lambda_1^{(3)}, d_*^{(2)},
\Lambda_1^{(2)},{\bar \Lambda}_1^{(2)}$ in (\ref{X-30-1}) are independent. ${\bar \Lambda}_1^{(2)}$
has the same distribution as $\Lambda_1^{(2)}$ .

We are ready to prove Theorem \ref{T4+}. The convergence  $C_G(k)\to C_*(k)$ as $n,m\to+\infty$ is shown in Theorem 2 of  \cite{BloznelisPetuchovas}.
Here we only prove (\ref{X-24-1}).
For $r=0,1,2,3$ we have, by Lemma \ref{Lambda}, 
\begin{equation}\label{X-30-1+}
\PP(\Lambda_0^{(r)}=t)
\sim 
c_0(r)
\,
t^{-(\gamma-r)},
\qquad
\PP(\Lambda_1^{(r)}=t)
\sim
c_1(r)
\,
t^{-(\alpha-r)}.
\end{equation}
 Furthermore, for $r=1,2$ we have, by Theorem \ref{T1+},
 \begin{equation}\label{X-30-2+}
 \PP(d_*^{(r)}=t)
 \sim 
 c_2(r,\alpha,\gamma)
 \,
 t^{-(\alpha-1)\wedge(\gamma-r)}.
 \end{equation}
We note that explicit expressions of  $c_0(r), c_1(r),
c_2(r,\alpha,\gamma)$ in terms  $a,b,\beta, a_i,b_i$ are easy to obtain, but we do not write down them here.
It follows from (\ref{X-30-1+}), (\ref{X-30-2+}) that 
\begin{eqnarray}\nonumber
\nonumber
p_1(t)
&
\sim 
&
{\mathbb I}_{\{\alpha\ge \gamma\}}
\PP(d_*^{(2)}=t)
+
{\mathbb I}_{\{\alpha\le \gamma\}}
\Bigl(
\PP(\Lambda_1^{(2)}=t)+\PP(\Lambda_2^{(2)}=t)
\Bigr)
\\
\nonumber
&
\sim
&
{\mathbb I}_{\{\alpha\ge \gamma\}}
c_2(2,\alpha,\gamma)t^{2-\gamma}
+
{\mathbb I}_{\{\alpha\le \gamma\}}
2c_1(2)t^{2-\alpha},
\\
\nonumber
p_2(t)
&
\sim 
&
{\mathbb I}_{\{\alpha\ge \gamma+2\}}
\PP(d_*^{(1)}=t)
+
{\mathbb I}_{\{\alpha\le \gamma+2\}}
\PP(\Lambda_1^{(3)}=t)
\\
\nonumber
&
\sim
&
{\mathbb I}_{\{\alpha\ge \gamma+2\}}
c_2(1,\alpha,\gamma)t^{1-\gamma}
+
{\mathbb I}_{\{\alpha\le \gamma+2\}}
c_1(3)t^{3-\alpha}.
\end{eqnarray}
 Combining these relations we conclude that $p_1(t)/p_2(t)$ scales as $t^{\varkappa}$, where
 $\varkappa=-1$ for $\alpha\le \gamma$, $\varkappa=\alpha-\gamma-1$ for $\gamma<\alpha<\gamma+2$, and $\varkappa=1$ for $\gamma+2\le \alpha$.
 Now (\ref{X-24-1}) follows from (\ref{X-28-2+}).
 $\qedsymbol$


\section{Auxiliary results}

 In Theorem \ref{Borovkovo} we collect several results 
from
Theorems 2.2.1, 2.2.3, 3.1.1, 3.1.6, 4.7.6  of \cite{Borovkov2008}, see also \cite{Mogulskii2008}. By $c_{X_1}(r_1,\dots, r_k)$ we denote a positive constant depending on the distribution of $X_1$ and numbers $r_1,\dots, r_k$.
 We observe that (\ref{V-04-11}) implies that for some $c_*>0$ 
 \begin{equation}\label{2019-01-09+1}
 \PP(X_1\ge t)\le c_* t^{1-\alpha}L_1(t),
\qquad
{\text{for}}
\qquad
t\ge 1.
\end{equation}

\begin{tm}\label{Borovkovo}
Let $\alpha, r\ge 1$.  Assume that (\ref{2019-01-09+1}) holds, where $L_1$ is slowly varying at infinity.


(i) Let $1<\alpha<2$. There exists  
$c=c_{X_1}(r)$ such  that for any $n\ge 1$ and $x\ge y>0$ satisfying
$x/y\le r$ we have
\begin{equation}\label{Borovkov1}
 \PP(S_n\ge x,\, M_n<y)
 \le 
 c\bigl(ny^{1-\alpha}L_1(y)\bigr)^{x/y}.
\end{equation}

(ii)
Let $\alpha=2$. For any $\beta, \eta>0$
 there exists 
$c=c_{X_1}(r, \beta,\eta)$ such  that (\ref{Borovkov1}) holds 
for any 
$n\ge 1$ and $x\ge y>0$ satisfying   
$x/y\le r$ and 
\begin{equation}\label{Borovkov1+}
n\PP(X_1\ge y)\bigl(L_\Delta(y)\bigr)^{1+\beta} 
\le \eta,
\qquad
{\text{where}}
\qquad
L_\Delta(y)=
\frac{\int_{0}^y\PP(X_1\ge u)du}{y\PP(X_1\ge y)}
\end{equation}
is a slowly varying function.

(iii)
Let $2\le \alpha<3$. Assume in addition that (\ref{2018-12-19-2})
holds and 
$\EE|X_1|<\infty$, $\EE X_1=0$. 
For $\alpha=2$ we also assume that 
$\PP(X_1<-t)=O(t^{-\beta})$ for some $\beta>1$.
Then 
there exists
$c=c_{X_1}(r)$ such  that for any $n\ge 1$ and $x\ge y>0$ satisfying
$x/y\le r$  inequality (\ref{Borovkov1}) holds.

(iv)
Let $\alpha=3$. Assume, in addition, that (\ref{2018-12-19-2})
holds and $\EE X_1=0$. Denote for $u>0$
\begin{eqnarray}
\nonumber
&&
V(u)
=
\,
\begin{cases}
u^{-2},
\qquad
\qquad
\qquad
\quad
\
\qquad
{\text{for}}
\quad
\int_0^{+\infty} t^{-1}L_1(t)dt<\infty,
\\
u^{-2}\int_{0}^{u}s^{-1}L_1(s)ds,
\quad
\qquad
{\text{for}}
\quad
\int_0^{+\infty} t^{-1}L_1(t)dt=\infty,
\end{cases}
\\
\nonumber
&&
W(u)
=
\begin{cases}
u^{-2},
\qquad
\qquad
\qquad\qquad
\qquad
\qquad
\quad
{\text{for}}
\quad
\int_0^{+\infty} t^{-1}L_1(t)dt<\infty,
\\
u^{-2}L_1(u)
\int_{0}^{u}
\bigl(\int_{s}^{\infty}t^{-2}L_1(t)dt\bigr)ds,
\quad
\
{\text{for}}
\quad
\int_0^{+\infty} t^{-1}L_1(t)dt=\infty.
\end{cases}
\end{eqnarray}
For any $\eta>0$ there exists
$c=c_{X_1}(r,\eta)$ such  that (\ref{Borovkov1}) holds for each
$n\ge 1$ and $x\ge y>0$ satisfying  
$x/y\le r$ 
and 
\begin{equation}\label{2019-01-10+2}
nV\bigl(y/|\ln \Pi(x)|\bigr)+nW\bigl(y/|\ln \Pi(x)|\bigr)<\eta,
\qquad
{\text{where}}
\qquad 
\Pi(x)=nc_*x^{-2}L_1(x).
\end{equation}

\noindent
(v)
Let $\alpha>3$. Assume that  (\ref{V-04-11}), (\ref{2019-01-06-1}) hold. Denote $\sigma^2=\Var X_1$.
We have 
uniformly~in~${t\ge\sqrt{n}}$ 
\begin{equation}\label{Borovkov8}
 \PP(S_n-\lfloor n\mu\rfloor =t)
\sim
\frac{1}{\sigma \sqrt{2\pi n}}
e^{-\frac{t^2}{2n\sigma^2}}
+
n(\alpha-1) t^{-1}\PP(X_1-\mu>t).
\end{equation}
\end{tm}
\noindent
Here $a_n(t)\sim b_n(t)$ uniformly in 
$t\ge \sqrt{n}$ 
means 
that
$\lim_{n\to\infty}\sup_{t\ge \sqrt{n}}a_n(t)/b_n(t)=1$.


\section{Appendix}
\begin{lemama}
\label{alpha=3+} 
Let $a,b\ge 0$  such that 
$a+b>0$.  Let $X_1,X_2,\dots$ be non-negative integer valued iid random variables such that
\begin{equation}\label{2019-02-16}
\PP(X_1=t)= (a+o(1))t^{-3},
\qquad
\PP(X_1=-t)= (b+o(1))t^{-3}
\qquad
{\text{as}}
\quad
t\to+\infty.
\end{equation}
Let  
$\{\eta_s,\, s=0,\pm 1, \pm 2,\dots\}$ be the sequence defined by 
\begin{displaymath}
\PP(X_1=t)=(a+\eta_t)t^{-3}
\qquad
{\text{and}}
\qquad
\PP(X_1=-t)=(b+\eta_{-t})t^{-3},
\qquad 
t=0,1,2\dots.
\end{displaymath}
Denote $\mu=\EE X_1$, 
$b_n=\sqrt{0.5 (a+b) n\ln n}$ and $h(k)=\sum_{1\le j\le k}(|\eta_j|+|\eta_{-j}|)/j$.
Let $\varphi(s)=(2\pi)^{-1/2}e^{-s^2/2}$ denote the standard normal density. There exist numbers $c, c_1>0$ independent of $t$ and $n$ such that for each $k=0,1,2,\dots$ and each $n=1,2,\dots$ we have
\begin{eqnarray}\label{alpha=3+1}
&&
\Bigl|\PP\bigl(X_1+\cdots+X_n=k\bigr)
-
\varphi\bigl(b_n^{-1}(k-n\mu)\bigr)
\Bigr|
\le 
c \min_{1<A<\ln^2 n}T(A),
\\
\nonumber
&&
T(A)
:=
A^5n^{-1}
+
A^3
\bigl(
h(\lfloor b_n\rfloor)+\ln\ln n
\bigr)
\ln^{-1}n
+
e^{-c_1A}.
\end{eqnarray}
\end{lemama}

\begin{rem}\label{X-11-1}
For 
$|\eta_s|
=O\bigl(  
(\ln\ln |s|)^{-1}(\ln\ln\ln |s|)^{-4}
\bigr)$ as $|s|\to+\infty$,
Lemma \ref{alpha=3+}  implies
\begin{equation}\label{alpha=3+++}
\Bigl|\PP\bigl(X_1+\cdots+X_n=k\bigr)
-
\varphi\bigl(b_n^{-1}(k-n\mu)\bigr)
\Bigr|=o\bigl(1/\ln\ln n\bigr).
\end{equation}
\end{rem}
Indeed,  we have for large $k$ that 
$
h(k)
\le
 \sum_{|j|\le k}|\eta_j|
\le 
c'(\ln k)/\bigl((\ln\ln k)(\ln\ln\ln k)^4\bigr)$. Now for $A=A_n=(\ln\ln\ln n)^{5/4}$ we obtain $T(A_n)=o\bigl(1/\ln\ln n\bigr)$.

	{\it Proof of Lemma \ref{alpha=3+}}.
The proof goes along the lines of the proof of Theorem 4.2.1 
of \cite{IbragimovLinnik}. We begin with introducing some notation.
 Denote $\Delta=\Delta_{n,k}$ the 
quantity on the left of (\ref{alpha=3+1}). 
Denote $f(t)=\EE e^{itX_1}$ and $\phi(t)=\EE e^{itY}$ the Fourier-Stieltjes 
transforms of the probability 
distributions of $X_1$ and $Y=b_n^{-1}(X_1-\mu)$ with $i$ standing for
the imaginary unit. We denote $D_t(y)=e^{ity}-1-ity$ and use the inequalities
$|D_t(y)|\le 2|ty|$ and $|D_t(y)+(ty)^2/2|\le|ty|^3/6$ for real numbers $t$ 
and $y$. 
Note that (\ref{2019-02-16}) implies $\eta(s)=o(1)$ as $|s|\to+\infty$.

Let us show (\ref{alpha=3+1}). Given $1<A<\pi b_n$ we put $\varepsilon=A^{-1}$.
We have 
(formula (4.2.5) of \cite{IbragimovLinnik}))
\begin{displaymath}
\Delta
\le
 I_1+I_2+I_3,
\end{displaymath}
where
\begin{equation}\nonumber
I_1
=
\int_{-A}^A
\bigl|
\phi^n(t)-e^{-t^2/2}\bigr|dt,
\qquad
I_2=\int_{A\le |t|\le \pi b_n}
\Bigl|f^n\Bigl(\frac{t}{b_n}\Bigr)\Bigr|dt,
\qquad
I_3=\int_{|t|\ge A}
e^{-t^2/2}dt.
\end{equation}
Furthermore, for any $0<\delta<2$ there 
exist $n_0>0$,  $c_{\delta}>0$ and 
$\varepsilon_{\delta}\in (0,1)$  such that 
for $n>n_0$ and 
$|t|\le \varepsilon_{\delta} b_n$ we have 
$\bigl|f^n(b_n^{-1}t)\bigr|\le e^{-
c_{\delta}|t|^{\delta}}$ 
(formula (4.2.7) of \cite{IbragimovLinnik}). 
We choose $\delta=1$. For 
$\varepsilon_1 b_n\le |t|\le \pi$ we have 
$\bigl|f^n(b_n^{-1}t)\bigr|\le e^{-c_*n}$, 
for some $c_*>0$ independent of $n$ (formula 
(4.2.9) of \cite{IbragimovLinnik}). 
These upper bounds for $\bigl|f^n(b_n^{-1}t)\bigr|$  imply the bound
\begin{equation}\label{alpha=3+2}
I_2\le 2c_1^{-1}e^{-c_1A}+2\pi b_ne^{-c_*n}.
\end{equation}
Next, we estimate 
$I_3\le 2A^{-1}e^{-A^2/2}$ using the  inequality
$\PP(W>A)\le A^{-1}e^{-A^2/2}(2\pi)^{-1/2}$ 
for the standard Gaussian random variable 
$W$, see Section 7.1 of \cite{Feller-I-1968}. 
Finally, we show that 
\begin{equation}\label{alpha=3+3}
I_1
\le 
c'A^5n^{-1}+c'\Bigl(
A^2\varepsilon^{-1}+A^3\bigl(|h(\lfloor b_n\rfloor)|+|\ln \varepsilon|+\ln\ln n\bigr)
+
A^4\varepsilon
\Bigr)\ln ^{-1}n.
\end{equation}
The bounds for $I_1,I_2,I_3$ above imply (\ref{alpha=3+1}). 
It remains to prove (\ref{alpha=3+3}).
The identity 
\begin{displaymath}
\phi^n(t)
-
e^{-t^2/2}
=
\bigl(\phi(t)-e^{-t^2/(2n)}\bigr)
\sum_{j=1}^{n}\phi^{n-j}(t)e^{-(j-1)t^2/(2n)}
\end{displaymath}
implies 
$\bigl|\phi^n(t)-e^{-t^2/2}\bigr|\le n\bigl|\phi(t)-e^{-t^2/(2n)}\bigr|=:n\Delta^*$. 
In order to estimate $\Delta^*$ we expand $\phi(t)$ and $e^{-t^2/(2n)}$ in powers of 
$t$. 
Note that $\EE Y=0$ implies $\phi(t)-1
=
\EE D_t(Y)$. We split
\begin{equation}\label{IX-29-3}
\EE D_t(Y)
=
\EE D_t(Y){\mathbb I}_{\{|Y|<\varepsilon \}}
+
\EE D_t(Y){\mathbb I}_{\{|Y|\ge \varepsilon \}}
=:
J_1+J_2.
\end{equation}
Using $|D_t(y)|\le 2|ty|$ we obtain 
$J_2
\le
 2|t|\EE|Y|{\mathbb I}_{\{|Y|\ge \varepsilon \}}$. A simple calculation shows that
 $\EE|Y|{\mathbb I}_{\{|Y|\ge \varepsilon \}}
 \le c'\varepsilon^{-1}b_n^{-2}$. Hence $J_2\le c'|t|\varepsilon^{-1}b_n^{-2}$. 
 Next, using $|D_t(y)+(ty)^2/2|\le |ty|^3/6$ we obtain
\begin{equation}\label{IX-29-2}
J_1=-2^{-1}t^2\EE Y^2{\mathbb I}_{\{|Y|<\varepsilon \}}
+
6^{-1}(it)^3R,
\qquad
|R|\le \EE |Y|^3{\mathbb I}_{\{|Y|<\varepsilon \}}.
\end{equation}
A calculation shows that $|R|\le c'\varepsilon b_n^{-2}$.  Furthermore, we have
\begin{eqnarray}\label{IX-29-1}
b_n^2\EE Y^2{\mathbb I}_{\{|Y|<\varepsilon\}}
&
=
&
\sum_{j:\, |j-\mu|\le\varepsilon b_n}
(j-\mu)^2\PP(X_1=j)
=
\sum_{0<|j|\le \varepsilon b_n}\frac{a+b}{j}+r
\\
\nonumber
&
=
&(a+b)\ln \bigl(\lfloor \varepsilon b_n\rfloor\bigr)+
r',
\end{eqnarray}
where $r,r'$ denote the remainders. We have
$|r'|\le h\bigl(\lfloor \varepsilon b_n\rfloor\bigr)+c'$.
Using the inequalities
\begin{eqnarray}
\nonumber
&&
\bigl|
(a+b)b_n^{-2}
\ln \lfloor \varepsilon b_n\rfloor -n^{-1}
\bigr|
\le 
c'\bigl(|\ln \varepsilon|+\ln\ln n\bigr)/(n\ln n),
\\
\nonumber
&&
\bigl|
h\bigl(\lfloor(\varepsilon b_n\rfloor\bigr)-h\bigl(\lfloor(b_n\rfloor\bigr)
\bigr|
\le 
c'\sum_
{\varepsilon b_n\le |j|\le b_n}
j^{-1}
\le c'|\ln \varepsilon|
\end{eqnarray}
we approximate  $(a+b)
\ln \lfloor \varepsilon b_n\rfloor$ by $b_n^2/n$  and $h\bigl(\lfloor(\varepsilon b_n\rfloor\bigr)$ by $h\bigl(\lfloor(b_n\rfloor\bigr)$ in 
(\ref{IX-29-1}).  From (\ref{IX-29-3}),  
(\ref{IX-29-2}),  (\ref{IX-29-1}) we obtain 
the expansion 
\begin{displaymath}
n\Bigl|\phi(t)-1+\frac{t^2}{2n}\Bigr|
\le \frac{c'}{\ln n}R^*,
\qquad
R^*
=
|t|\varepsilon^{-1}+t^2\bigl(\ln\ln n+|\ln \varepsilon|+
h(\lfloor b_n\rfloor)\bigr)
+
|t|^3\varepsilon.
\end{displaymath}
We compare it with the  expansion 
$n\bigl|e^{-t^2/(2n)}-1+t^2/(2n)\bigr|\le t^4/(4n)$ and conclude that 
$n\Delta^*\le c'R^*\ln^{-1} n+t^4/(4n)$. This inequality implies
(\ref{alpha=3+3}).
$\qedsymbol$

\begin{lemama}
\label{Lambda}
Let $\alpha>2$, $a,b>0$ and let  $k$ be a positive integer. Let $Z$ be a non-negative integer valued random variable such that $\PP(Z=t)\sim at^{-\alpha}$ as $t\to+\infty$. Then
\begin{equation}\label{X-27-1}
\EE \Bigl(\frac{e^{-bZ}(bZ)^t}{t!}\Bigr)
\sim 
a \, b^{\alpha-1}t^{-\alpha}
\qquad
{\text{as}}
\quad 
t\to+\infty.
\end{equation} 
\end{lemama}
	{\it Proof of Lemma \ref{Lambda}}.
 Denote ${\tilde Z}=bZ$, $f_t(\lambda)=e^{-\lambda}\lambda^t(t!)^{-1}$ and $u_t=t^{1/2}\ln t$.
We split
\begin{displaymath}
\EE f_t({\tilde Z})
=
\EE f_t({\tilde Z}){\mathbb I}_{\{|{\tilde Z}-t|\le u_t\}}
+
\EE f_t({\tilde Z}){\mathbb I}_{\{{\tilde Z}<t-u_t\}}
+
\EE f_t({\tilde Z}){\mathbb I}_{\{{\tilde Z}>t+u_t\}}
=: I_1+I_2+I_3
\end{displaymath}
and show that $I_1\sim ab^{\alpha-1}t^{-\alpha}$ and
$I_j=o(t^{-\alpha})$ for $j=2,3$.

Let $\eta_1,\eta_2,\dots$ be iid Poisson random 
variables with mean $b$. By Theorem 6 chpt. 7 of
\cite{Petrov}, relation (\ref{VIII-23-3}) holds for  
the sum ${\tilde S}_n=\eta_1+\dots+\eta_n$
with $b_n=\sqrt {bn}$, $a_n=bn$ and  
$\tau_n=O(n^{-1/2})$.
Now from Lemma \ref{lema2++} (which applies 
to the sum ${\tilde S}_n$ as well) we obtain
\begin{displaymath}
\sum_{n:\,|bn-t|\le u_t}f_t(bn)
=
\sum_{n:\,|bn-t|\le u_t}\PP({\tilde S}_n=t)
\to b^{-1}.
\end{displaymath}
This relation implies $I_1\sim ab^{\alpha-1}t^{-\alpha}$.
The remaining  bounds $I_j=o(t^{-\alpha})$, $j=2,3$ are 
easy. Using the fact that 
$\lambda\to f_t(\lambda)$ increases (decreases) for $\lambda<t$ (for 
$\lambda>t$) we obtain for large 
$t$
\begin{eqnarray}
\nonumber
I_2
&
\le 
&
f_t(t-u_t)
\le 
e^{u_t}\Bigl(1-\frac{u_t}{t}\Bigr)^t
\le e^{-0.5\ln^2t},
\\
\nonumber
I_3
&
\le 
& 
f_t(t+u_t)
\le 
e^{-u_t}\Bigl(1+\frac{u_t}{t}\Bigr)^t
\le
e^{-0.1\ln^2t}.
\end{eqnarray}
Here we  applied $t!\ge (t/e)^t$ and then evaluated
$
\Bigl(1\pm \frac{u_t}{t}\Bigr)^t
=
e^{t\ln(1\pm u_tt^{-1})}
$
using a two term expansion of $\ln(1\pm u_tt^{-1})$ in powers of $u_tt^{-1}$.
$\qedsymbol$

{\it Acknowledgment.}
I thank anonymous referees for comments and remarks. Especially, for
suggesting the use of the discrete 
renewal theorem which makes the result of Theorem \ref{T1+} stronger and the proof 
 simpler. I am also grateful for pointing out the  papers
\cite{Doney1989}, \cite{EmbrechtsMaejimaOmey1984},
\cite{Hilberdink1996}.

\end{document}